\documentclass[11pt]{article}
\usepackage{amsfonts}
\usepackage{amssymb}
\usepackage{amsmath}
\def\sqr#1#2{{\vcenter{\vbox{\hrule height.#2pt
              \hbox{\vrule width.#2pt height#1pt \kern#1pt \vrule
width.#2pt}
              \hrule height.#2pt}}}}
\def\signed #1{{\unskip\nobreak\hfil\penalty50
              \hskip2em\hbox{}\nobreak\hfil#1
              \parfillskip=0pt \finalhyphendemerits=0 \par}}
\def\endpf{\signed {$\sqr69$}}
\def\dbC{{\mathbb{C}}}

\def\dbE{{\mathbb{E}}}
\def\dbF{{\mathbb{F}}}

\def\dbH{{\mathbb{H}}}

\def\dbN{{\mathbb{N}}}

\def\dbP{{\mathbb{P}}}

\def\dbR{{\mathbb{R}}}

\def\b{\beta}

\def\d{\delta}
\def\e{\varepsilon}

\def\si{\sigma}

\def\f{\varphi}

\def\o{\omega}

\def\3n{\negthinspace \negthinspace \negthinspace }
\def\2n{\negthinspace \negthinspace }
\def\1n{\negthinspace }
\def\ns{\noalign{\smallskip} }

\def\ds{\displaystyle}
\def\G{\Gamma}
\def\D{\Delta}

\def\Si{\Sigma}

\def\O{\Omega}
\def\cA{{\cal A}}

\def\cD{{\cal D}}

\def\cF{{\cal F}}
\def\cG{{\cal G}}
\def\cH{{\cal H}}

\def\cJ{{\cal J}}

\def\cL{{\cal L}}

\def\cQ{{\cal Q}}

\def\cT{{\cal T}}
\def\cU{{\cal U}}

\def\cY{{\cal Y}}

\def\mE{{\mathbb{E}}}

\def\no{\noindent}

\def\ms{\medskip}
\def\bs{\bigskip}
\def\q{\quad}
\def\qq{\qquad}
\def\hb{\hbox}

\def\lan{\mathop{\langle}}
\def\ran{\mathop{\rangle}}

\def\pa{\partial}

\def\wt{\widetilde}
\def\cd{\cdot}
\def\cds{\cdots}

\def\dim{\hbox{\rm dim$\,$}}

\def\ae{\hbox{\rm a.e.{ }}}

\def\span{\hbox{\rm span$\,$}}

\def\Re{{\mathop{\rm Re}\,}}

\def\({\Big (}
\def\){\Big )}
\def\[{\Big[}
\def\]{\Big]}

\def\={\buildrel \triangle \over =}

\def\resp{{\it resp. }}

\def\be{\begin{equation}}
\def\bel{\begin{equation}\label}
\def\ee{\end{equation}}
\def\bea{\begin{eqnarray}}
\def\eea{\end{eqnarray}}
\def\bt{\begin{theorem}}
\def\et{\end{theorem}}
\def\bc{\begin{corollary}}
\def\ec{\end{corollary}}
\def\bl{\begin{lemma}}
\def\el{\end{lemma}}
\def\bp{\begin{proposition}}
\def\ep{\end{proposition}}
\def\br{\begin{remark}}
\def\er{\end{remark}}
\def\ba{\begin{array}}
\def\ea{\end{array}}
\def\bd{\begin{definition}}
\def\ed{\end{definition}}

\newtheorem{lemma}{Lemma}[section]
\newtheorem{remark}{Remark}[section]

\newtheorem{theorem}{Theorem}[section]
\newtheorem{corollary}{Corollary}[section]

\newtheorem{definition}{Definition}[section]
\newtheorem{proposition}{Proposition}[section]

\makeatletter
   
   \@addtoreset{equation}{section}
\makeatother

\begin{document}

\title{\bf Transposition Method for Backward Stochastic Evolution Equations Revisited,
and Its Application}

\author{Qi L\"{u}\thanks{School of Mathematics, Sichuan University, Chengdu 610064, China. {\small\it e-mail:} {\small\tt
luqi59@163.com}. \ms}
 ~~~
  and~~~
Xu Zhang\thanks{Yangtze Center of
Mathematics, Sichuan University, Chengdu
610064, China. {\small\it e-mail:}
{\small\tt zhang$\_$xu@scu.edu.cn}.}}
\date{}

\maketitle

\begin{abstract}
The main purpose of this paper is to
improve our transposition method to solve
both vector-valued and operator-valued
backward stochastic evolution equations
with a general filtration. As its
application, we obtain a general
Pontryagin-type maximum principle for
optimal controls of stochastic evolution
equations in infinite dimensions. In
particular, we drop the technical
assumption appeared in \cite[Theorem
9.1]{LZ1}.
\end{abstract}

\bs

\no{\bf 2010 Mathematics Subject
Classification}. 93E20

\bs

\no{\bf Key Words}. Backward stochastic
evolution equations, transposition
solution, optimal control, Pontryagin-type
maximum principle.


\section{Introduction}\label{s1}


Let $T>\tau\ge 0$, $d\in\dbN$, and
$(\O,\cF,\{\cF_t\}_{t\in [0,T]},\dbP)$ be a
complete filtered probability space
(satisfying the usual conditions), on which
a standard $d$-dimensional Brownian motion
$\{w(t)\}_{t\in [0,T]}$ is defined. Write
$\dbF=\{\cF_t\}_{t\in [0,T]}$, and let $X$
be a Banach space. For any $t\in[0,T]$ and
$r\in [1,\infty)$, denote by
$L_{\cF_t}^r(\O;X)$ the Banach space of all
$\cF_t$-measurable random variables
$\xi:\O\to X$ such that
$\mathbb{E}|\xi|_X^r < \infty$, with the
canonical  norm.  Also, denote by
$L^{r}_{\dbF}(\O;D([\tau,T];X))$ the vector
space of all $X$-valued c\`adl\`ag process
$\phi(\cd)$ such that $\mE\big(\sup_{t\in
[\tau,T)}|\phi(r)|_X\big)^r<\infty$. One
can show that
$L^{r}_{\dbF}(\O;D([\tau,T];X))$ is a
Banach space with the norm
 $$
|\phi(\cd)|_{L^{r}_{\dbF}(\O;D([\tau,T];X))}
\= \big[\mE (\sup_{t\in
[\tau,T)}|\phi(r)|_X)^r\big]^{\frac{1}{r}}.
 $$
Denote by $C_{\dbF}([\tau,T];L^{r}(\O;X))$
the Banach space of all $X$-valued
$\dbF$-adapted processes $\phi(\cdot)$ such
that $\phi(\cdot):[\tau,T] \to
L^{r}_{\cF_T}(\O;X)$ is continuous, with
the norm
$$
|\phi|_{C_{\dbF}([\tau,T];L^{r}(\O;X))}\=\sup_{t\in
[\tau,T)}\big(\mE
|\phi(r)|_X^r\big)^{\frac{1}{r}}.
$$
Similarly, one can define the Banach space
$D_{\dbF}([\tau,T];L^{r}(\O;X))$.

Fix $r_1,r_2,r_3,r_4\in[1,\infty]$. Put

 $$
 \ba{ll}
\ds
L^{r_1}_\dbF(\O;L^{r_2}(\tau,T;X))=\Big\{\f:(\tau,T)\times\O\to
X\bigm|\f(\cd)\hb{
is $\dbF$-adapted and }\dbE\(\int_\tau^T|\f(t)|_X^{r_2}dt\)^{\frac{r_1}{r_2}}<\infty\Big\},\\
\ns\ds
L^{r_2}_\dbF(\tau,T;L^{r_1}(\O;X))=\Big\{\f:(\tau,T)\times\O\to
X\bigm|\f(\cd)\hb{ is $\dbF$-adapted and
}\int_\tau^T\(\dbE|\f(t)|_X^{r_1}\)^{\frac{r_2}
{r_1}}dt<\infty\Big\}.
 \ea
 $$
Both $L^{r_1}_\dbF(\O;L^{r_2}(\tau,T;X))$
and $L^{r_2}_\dbF(\tau,T;L^{r_1}(\O;X))$
are Banach spaces with the canonical norms.
If $r_1=r_2$, we simply denote the above
spaces by $L^{r_1}_\dbF(\tau,T;X)$. Let $Y$
be another Banach space. Denote by $\cL(X,
Y)$ the (Banach) space of all bounded
linear operators from $X$ to $Y$, with the
usual operator norm (When $Y=X$, we simply
write $\cL(X)$ instead of $\cL(X, Y)$).
Further, we denote by
$\cL_{pd}\big(L^{r_1}_{\dbF}(\tau,T;L^{r_2}(\O;X)),\;L^{r_3}_{\dbF}(\tau,T;L^{r_4}(\O;Y))\big)$
(\resp
$\cL_{pd}\big(X,\;L^{r_3}_{\dbF}(\tau,T;L^{r_4}(\O;Y))\big)$)
the vector space of all bounded,
pointwisely defined linear operators $G$
from $L^{r_1}_{\dbF}(\tau,T;L^{r_2}(\O;X))$
(\resp $X$) to
$L^{r_3}_{\dbF}(\tau,T;L^{r_4}(\O;Y))$,
i.e., for $\ae (t,\omega)\in
(\tau,T)\times\Omega$, there exists an
$L(t,\o)\in\cL (X,Y)$ such that $\big(G
u(\cd)\big)(t,\o)=L (t,\o)u(t,\o), \;\,
\forall\; u(\cd)\in
L^{r_1}_{\dbF}(\tau,T;L^{r_2}(\O;X))$
(\resp $\big(G x\big)(t,\o)=L (t,\o)x, \;\,
\forall\; x\in X$). In a similar way, one
can define the spaces such as
$\cL_{pd}\big(L^{r_2}_{\cF_t}(\O;X),\;L^{r_3}_{\dbF}(\tau,T;L^{r_4}(\O;Y))\big)$
and
$\cL_{pd}\big(L^{r_2}_{\cF_t}(\O;X),L^{r_4}_{\cF_t}(\O;Y)\big)$
for $t\in [0,T]$, etc.

Let $H$ be a complex Hilbert space. Denote
by $H^d$ the Cartesian product
$\underbrace{H\times H\times \cdots\times
H}_{d\hb{ \tiny times}}$. Similarly, we
will use the notations $\cL(H)^d$,
$\cL_2(H)^d$ and so on, where  $\cL_2(H)$
stands for the (Hilbert) space of all
Hilbert-Schmidt operators on $H$.

 Let $A$ be an unbounded linear operator
(with domain $D(A)\subset H$), which
generates a $C_0$-semigroup
$\{S(t)\}_{t\geq 0}$ on $H$. Denote by
$A^*$ the dual operator of $A$. It is
well-known that $D(A)$ is a Hilbert space
with the usual graph norm, and $A^*$
generates a $C_0$-semigroup
$\{S^*(t)\}_{t\geq 0}$, which is the dual
$C_0$-semigroup of $\{S(t)\}_{t\geq 0}$.

First, we consider the following $H$-valued
backward stochastic evolution equation
(BSEE for short):
\begin{eqnarray}\label{bsystem1}
\left\{
\begin{array}{lll}
\ds dy(t) = -  A^* y(t) dt + f(t,y(t),Y(t))dt + Y(t) dw(t) &\mbox{ in }[\tau,T),\\
\ns\ds y(T) = y_T,
\end{array}
\right.
\end{eqnarray}
where $y_T \in L_{\cF_T}^{p}(\O;H)$ with
$p\in (1,\infty]$, and
$f(\cd,\cd,\cd):[\tau,T]\times H\times H^d
\to H$ satisfies, for some constant
$C_L>0$,
\begin{equation}\label{Lm1}
\left\{
\begin{array}{ll}\ds
f(\cd,0,0)\in
L^{1}_{\dbF}(\tau,T;L^{p}(\O;H)),\\\ns\ds
|f(t,x_1,y_1)-f(t,x_2,y_2)|_H\leq
C_L\big(|x_1-x_2|_H+|y_1-y_2|_{H^d} \big),\\
\ns\ds\hspace{3.5cm} \ae (t,\o)\in
[\tau,T]\times\O,\;\; \forall\;x_1,x_2\in
H,\,y_1,y_2\in H^d.
\end{array}
\right.
 \end{equation}
Here neither the usual natural filtration
condition nor the quasi-left continuity is
assumed for the filtration $\dbF$, and  the
unbounded operator $A$ is only assumed to
generate a general $C_0$-semigroup. Hence,
we cannot apply the existing results on
infinite dimensional BSEEs (e.g.
\cite{Al-H1,HP2,MY,MM}) to obtain the
well-posedness of the equation
\eqref{bsystem1}.

Next, we consider the following
$\cL(H)$-valued BSEE\footnote{Throughout
this paper, for any operator-valued process
(\resp random variable) $R$, we denote by
$R^*$ its pointwisely dual operator-valued
process (\resp random variable). For
example, if $R\in L^{r_1}_\dbF(\tau,T;
L^{r_2}(\O; \cL(H)))$, then $R^*\in
L^{r_1}_\dbF(\tau,T; L^{r_2}(\O; \cL(H)))$,
and $|R|_{L^{r_1}_\dbF(\tau,T; L^{r_2}(\O;
\cL(H)))}=|R^*|_{L^{r_1}_\dbF(\tau,T;
L^{r_2}(\O; \cL(H)))}$.}:
\begin{equation}\label{bsystem2}
\left\{\3n
\begin{array}{ll}
\ds dP  =  - (A^*  + J^* )P dt  -  P(A + J
)dt -K^*PKdt
 - (K^* Q +  Q K)dt
  +   Fdt  +  Q dw(t) &\mbox{ in } [\tau,T),\\
\ns\ds P(T) = P_T,
\end{array}
\right.
\end{equation}
where
\begin{equation}\label{JKFP}
\begin{array}{ll}\ds
J\in L^4_\dbF(\tau,T; L^\infty(\O;
\cL(H))),\q K\in L^4_\dbF(\tau,T;
L^\infty(\O; \cL(H)^d)), \\
\ns\ds F\in
L^1_\dbF(\tau,T;L^p(\O;\cL(H))),\q P_T\in
L^p_{\cF_T}(\O;\cL(H)).
\end{array}
\end{equation}
If $H=\dbR^m$ for some $m\in\dbN$, then
\eqref{bsystem2} is an $m\times m$
matrix-valued backward stochastic
differential equation (BSDE for short), and
hence, one can easily obtain its
well-posedness for this special case. On
the other hand, if $\dim H=\infty$, $F\in
L^1_\dbF(\tau,T;L^p(\O;\cL_2(H)))$ and
$P_T\in L^p_{\cF_T}(\O;\cL_2(H))$, then
\eqref{bsystem2} is a special case of
\eqref{bsystem1} (because $\cL_2(H)$ is a
Hilbert space), and therefore in this case
the well-posedness of \eqref{bsystem2}
follows from that of \eqref{bsystem1}.
However, the situation is completely
different when $\dim H=\infty$ if one does
not impose further assumptions on $F$ and
$P_T$. Indeed, in the infinite dimensional
setting, although $\cL(H)$ is still a
Banach space, it is neither reflexive nor
separable even if $H$ itself is separable.
Because of this, $\cL(H)$ is NOT a UMD
space (needless to say a Hilbert space),
and consequently, it is even a quite
difficult problem to define the stochastic
integral ``$\int_\tau^T Q dw(t)$" (appeared
in (\ref{bsystem2})) for an $\cL(H)$-valued
process $Q$. We refer to \cite{Du, LZ1} for
previous studies on the well-posedness of
\eqref{bsystem2}, by avoiding the
definition of ``$\int_\tau^T Q dw(t)$" in
one way or another.

Similar to the finite dimensional case
(\cite{Peng1}), both \eqref{bsystem1} and
\eqref{bsystem2} play crucial roles in
establishing the Pontryagin-type maximum
principle for optimal controls of general
infinite dimensional nonlinear stochastic
systems with control-dependent diffusion
terms and possibly nonconvex control
regions (\cite{Du, Fuhrman, LZ1, TL}).

The main purpose of this paper is to
improve our transposition method, developed
in \cite{LZ1}, to solve the equations
\eqref{bsystem1} and \eqref{bsystem2}.
Especially, we shall give some
well-posedness/regularity results for
solutions to these two equations such that
they can be conveniently used in the above
mentioned Pontryagin-type maximum
principle.  In the stochastic finite
dimensional setting, the transposition
method (for solving BSDEs) was introduced
in our paper \cite{LZ}, but one can find a
rudiment of this method at \cite[pp.
353--354]{YZ}.

We remark that, our method is also
motivated by the classical transposition
method to solve the non-homogeneous
boundary value problems for deterministic
partial differential equations (see
\cite{Lions1} for a systematic introduction
to this topic) and especially the boundary
controllability problem for hyperbolic
equations (\cite{Lions}).

For the readers' convenience, let us recall
below the main idea of the classical
transposition method to solve the following
deterministic wave equation with
non-homogeneous Dirichelt boundary
conditions:
 \bel{ch2-s1}
 \left\{\ba{ll}
 \displaystyle
 y_{tt}-\Delta y=0 & \mbox{ in } Q\= (0,T)\times G, \\
 \displaystyle
y= u& \mbox{ on }\Si\=(0,T)\times \G, \\
\displaystyle
y(0)=y_0, \q y_t(0)=y_1& \mbox{ in } G,\\
 \ea\right.
 \ee
where $G$ is a nonempty open bounded domain
in $\dbR^d$ with a $C^2$ boundary $\G$,
$(y_0,y_1)\in L^2(G)\times H^{-1}(G)$ and
$u\in L^2(\Si)$ are given, and $y$ is the
unknown.

When $u\equiv 0$, one can use the standard
semigroup theory to prove the
well-posedness of (\ref{ch2-s1}) in the
solution space $ C([0,T];L^2(G))\bigcap
C^1([0,T];H^{-1}(G))$.

When $u\not\equiv 0$, one needs to use the
transposition method because $y|_\Si= u$
does NOT make sense by the usual trace
theorem. For this purpose, for any $f\in
L^1(0,T;L^2(G))$ and $g\in L^1(0,T;
H^1_0(G))$, consider the following adjoint
equation of (\ref{ch2-s1}):
 \be\label{ch23-1111}
 \left\{ \ba{ll}
\ds \zeta_{tt}-\Delta\zeta=f+g_t,\q&\hbox{ in }Q,\\[2mm]\ms
\zeta=0,&\hbox{ on }\Si,\\\ms
 \zeta(T)=\zeta_t(T)=0,&\hbox{ in }G.
  \ea\right.
  \ee
This equation admits a unique solution
$\zeta (\in C([0,T];$ $H^1_0(G))\bigcap
C^1([0,T];L^2(G)))$, which enjoys a hidden
regularity $\frac{\partial\zeta}{\partial
\nu}\in L^2(\Si)$ (c.f. \cite{Lions}). Here
and henceforth, $\nu\equiv \nu(x)$ stands
for the unit outward normal vector of $G$
at $x\in\G$.

In order to give a reasonable definition
for the solution to (\ref{ch2-s1}) by the
transposition method, we consider first the
case when $y$ is sufficiently smooth.
Assume that $g\in C_0^\infty(0,T;
H^1_0(G))$, $y_1\in L^2(G)$, and that $y\in
H^2(Q)$ satisfies (\ref{ch2-s1}). Then,
multiplying the first equation in
(\ref{ch2-s1}) by $\zeta$, integrating it
in $Q$, and using integration by parts, we
find that
  \bel{fkyk}\ba{ll}
\ds \int_Q fydxdt-\int_Q gy_tdxdt
=\int_G\zeta(0)y_1dx-\int_G\zeta_t(0)
y_0dx-\int_{\Si}\frac{\pa\zeta}{\pa
\nu}ud\Si.
 \ea
 \ee
Note that (\ref{fkyk}) still makes sense
even if the regularity of $y$ is relaxed as
$y\in C([0,T];L^2(G))$ $\bigcap
C^1([0,T];H^{-1}(G))$. This leads to the
following notion:

\bd\label{1d1}
We call $y\in C([0,T];L^2(G))\bigcap
C^1([0,T];H^{-1}(G))$ a transposition
solution to (\ref{ch2-s1}), if $y(0)=y_0$,
$y_t(0)=y_1$, and for any $f\in
L^1(0,T;L^2(G))$ and $g\in L^1(0,T;
H^1_0(G))$, it holds that
 $$\ba{ll}
 \ds \int_Q fydxdt-\int_0^T{\lan g,y_t\ran}_{H_0^1(G),H^{-1}(G)}dt ={\lan\zeta(0),y_1\ran}_{H_0^1(G),H^{-1}(G)}+\int_\O\zeta_t(0)
y_0dx-\int_{\Si}\frac{\pa\zeta}{\pa
\nu}ud\Si,
 \ea
 $$
where $\zeta$ is the unique solution to
(\ref{ch23-1111}). \ed

One can show the well-posedness of
(\ref{ch2-s1}) in the sense of Definition
\ref{1d1}, by means of the transposition
method (\cite{Lions}). Clearly, the point
of this method is to interpret the solution
to a forward wave equation with
non-homogeneous Dirichlet boundary
conditions in terms of another backward
wave equation with non-homogeneous source
terms. Of course, in the deterministic
setting, since the wave equation is
time-reversible, there exists no essential
difference between the forward problem and
the backward one. Nevertheless, this
reminds us to interpret BSDEs/BSEEs in
terms of forward stochastic
differential/evolution equations, as we
have done in \cite{LZ, LZ1}. Clearly, the
transposition method is a variant of the
standard duality method, and in some sense
it provides a way to see something which is
not easy to be detected directly.

The rest of this paper is organized as
follows. Section \ref{sec well1}  is
addressed to the well-posedness of the
equation \eqref{bsystem1}. Sections
\ref{sec well2} and \ref{sec regular} are
devoted to the well-posedness of the
equation \eqref{bsystem2} and a regularity
property for its solutions, respectively.
Finally, in Section \ref{max}, we show a
stochastic Pontryagin-type maximum
principle for controlled stochastic
evolution equations in infinite dimensions.

\section{Well-posedness of the vector-valued BSEEs}
\label{sec well1}

In this section, we discuss the
well-posedness of the equation
\eqref{bsystem1} in the transposition
sense.

Consider the following (forward) stochastic
evolution equation:
\begin{equation}\label{fsystem1}
\left\{
\begin{array}{lll}\ds
dz = (Az + \psi_1)ds +  \psi_2 dw(s) &\mbox{ in }(t,T],\\
\ns\ds z(t)=\eta,
\end{array}
\right.
\end{equation}
where $t\in[\tau,T]$, $q=\frac{p}{p-1}$,
$\psi_1\in L^q_{\dbF}(\O;L^1(t,T;H))$,
$\psi_2\in L^2_{\dbF}(t,T;L^q(\O; H^d))$
and $\eta\in L^{q}_{\cF_t}(\O;H)$. Let us
recall that $z(\cd)\in
C_\dbF([t,T];L^q(\O;H))$ is a (mild)
solution to the equation \eqref{fsystem1}
if
$$
z(s) = S(s-t)\eta + \int_t^s
S(s-\si)\psi_1(\si)d\si + \int_t^s
S(s-\si)\psi_2(\si)dw(\si),\;\forall\, s\in
[t,T].
$$
We now introduce the following notion.

\begin{definition}\label{definition1}
We call $(y(\cdot), Y(\cdot)) \in
L^{p}_{\dbF}(\O;D([\tau,T];H)) \times
L^2_{\dbF}(\tau,T;L^p(\O;H^d))$  a
transposition solution to
 \eqref{bsystem1} if for any $t\in
[\tau,T]$, $\psi_1(\cdot)\in
L^q_{\dbF}(\O;L^1(t,T;H))$,
$\psi_2(\cdot)\in L^2_{\dbF}(t,T;L^q(\O;
H^d))$, $\eta\in L^q_{\cF_t}(\O;H)$ and the
corresponding solution $z\in
C_{\dbF}([t,T];L^q(\O;H))$ to
\eqref{fsystem1}, it holds that
\begin{equation}\label{eq def solzz}
\begin{array}{ll}\ds
\dbE \big\langle z(T),y_T\big\rangle_{H}
- \dbE\int_t^T \big\langle z(s),f(s,y(s),Y(s) )\big\rangle_Hds\\
\ns\ds = \dbE
\big\langle\eta,y(t)\big\rangle_H +
\dbE\int_t^T \big\langle
\psi_1(s),y(s)\big\rangle_H ds +
\dbE\int_t^T \big\langle
\psi_2(s),Y(s)\big\rangle_{H^d} ds.
\end{array}
\end{equation}
\end{definition}

In what follows, we will use $C$ to denote
a generic positive constant, which may be
different from one place to another. We
have the following result for the
well-posedness of the equation
\eqref{bsystem1}.

\begin{theorem}\label{vw-th1}
For any  $y_T \in L^p_{\cF_T}(\O;H)$,
$f(\cd,\cd,\cd):[\tau,T]\times H\times H^d
\to H$ satisfying \eqref{Lm1}, the equation
\eqref{bsystem1} admits one and only one
transposition solution $(y(\cdot),
Y(\cdot))\in L^{p}_{\dbF}(\O;D([\tau,T];H))
\times L^2_{\dbF}(\tau,T;L^p(\O;H^d))$.
Furthermore,
\begin{equation}\label{vw-th1-eq1}
\begin{array}{ll}\ds
|(y(\cdot), Y(\cdot))|_{
L^{p}_{\dbF}(\O;D([t,T];H)) \times L^2_{\dbF}(\tau,T;L^p(\O;H^d))}\\
\ns\ds\leq C\left[
 |f(\cd,0,0)|_{ L^1_{\dbF}(t,T;L^p(\O;H))} +|y_T|_{
L^p_{\cF_T}(\O;H)}\right], \q\forall\;t\in
[\tau,T].
\end{array}
\end{equation}
\end{theorem}

\br
In \cite[Theorem 3.1]{LZ1}, it was assumed
that $1<p\le 2$ and $n=1$.  Also, in
Theorem \ref{vw-th1} the solution space of
the first unknown $y$ in \eqref{bsystem1}
is $L^{p}_{\dbF}(\O;D([\tau,T];H))$; while
in \cite{LZ1} this space is
$D_{\dbF}([\tau,T];L^p(\O;H))$. It is easy
to see that
$$
L^{p}_{\dbF}(\O;D([\tau,T];H))\hookrightarrow
D_{\dbF}([\tau,T];L^p(\O;H)),
$$
algebraically and topologically. Hence,
compared to \cite[Theorem 3.1]{LZ1},
Theorem \ref{vw-th1} improves a little the
regularity of solutions to
\eqref{bsystem1}.
\er

Before proving Theorem \ref{vw-th1}, we
first recall the following Riesz-type
Representation Theorem (See \cite[Corollary
2.3 and Remark 2.4]{LYZ}).

\begin{lemma}\label{lemma1}
Fix $t_1$ and $t_2$ satisfying $0\leq t_2 <
t_1 \leq T$.  Assume that $\cY$ is a
reflexive Banach space. Then, for any
$r,s\in [1,\infty)$, it holds that

$$
\left(L^r_\dbF(t_2,t_1;L^s(\O;\cY))\right)^*=L^{r'}_\dbF(t_2,t_1;L^{s'}(\O;\cY^*)),
$$
and
$$
\left(L^s_\dbF(\O;L^r(t_2,t_1;\cY))\right)^*=L^{s'}_\dbF(\O;L^{r'}(t_2,t_1;\cY^*)),
$$
 where
$$
s'=\left\{
\begin{array}{ll}\ds
s/(s-1), &\mbox{ if } s\not=1,\\
\ns\ds \infty &\mbox{ if } s\not=1;
\end{array}
\right. \qq r'=\left\{
\begin{array}{ll}\ds
r/(r-1), &\mbox{ if } r\not=1,\\
\ns\ds \infty &\mbox{ if } r\not=1.
\end{array}
\right.
$$

\end{lemma}

{\it Proof of Theorem \ref{vw-th1}}\,: It
suffices  to consider a particular case for
\eqref{bsystem1}, i.e. the case that
$f(\cd,\cd,\cd)$ is independent of the
second and third arguments. More precisely,
we consider the following equation:
\begin{eqnarray}\label{bsystem1-1}
\left\{
\begin{array}{lll}
\ds dy (t)= -  A^* y(t) dt + f(t)dt + Y (t)dw (t)&\mbox{ in }[\tau,T),\\
\ns\ds y(T) = y_T,
\end{array}
\right.
\end{eqnarray}
where $y_T \in L^p_{\cF_T}(\O;H)$ and
$f(\cdot)\in L^1_{\dbF}(\tau,T;
L^p(\O;H))$. The general case follows from
the well-posedness for (\ref{bsystem1-1})
and the standard fixed point technique.

We divide the proof into several steps.
Since the proof is very similar to that of
\cite[Theorem 3.1]{LZ1}, we give below only
a sketch.

\ms

{\bf Step 1.} For any $t\in [\tau,T]$, we
define a linear functional $\ell$
(depending on $t$) on the Banach space
$L^q_{\dbF}(\O;L^1(t,T;H))\times
L^2_{\dbF}(t,T;L^q(\O;H^d))\times
L^q_{\cF_t}(\O;H)$ as follows:
\begin{equation}\label{vw-th1-eq1.1}
\begin{array}{ll}
\ds\ell\big(\psi_1(\cdot),
\psi_2(\cdot),\eta\big) =
\mathbb{E}\big\langle z(T),y_T\big\rangle_H
- \mathbb{E}\int_t^T \big\langle
z(s),f(s)\big\rangle_H ds,\\\ns\ds
\qq\forall\; \big(\psi_1(\cdot),
\psi_2(\cdot),\eta\big)\in
L^q_{\dbF}(\O;L^1(t,T;H))\times
L^2_{\dbF}(t,T;L^q(\O;H^d))\times
L^q_{\cF_t}(\O;H),
\end{array}
\end{equation}
where $z(\cdot)\in
C_{\dbF}([t,T];L^{q}(\O;H))$ solves the
equation \eqref{fsystem1}. It is an easy
matter to show that $\ell$ is a bounded
linear functional on
$L^q_{\dbF}(\O;L^1(t,T;H))\times
L^2_{\dbF}(t,T;L^q(\O;H^d))\times
L^q_{\cF_t}(\O;H)$. By Lemma \ref{lemma1},
there exists a triple
$$
\big(y^t(\cdot), Y^t(\cdot), \xi^t\big)\in
L^p_{\dbF}(\O;L^\infty(t,T;H))\times
L^2_{\dbF}(t,T;L^p(\O;H^d))\times
L^p_{\cF_t}(\O;H)
$$
such that
\begin{equation}\label{vw-th1-eq2}
\begin{array}{ll}\ds
\mathbb{E}\big\langle z(T),y_T\big\rangle_H
-
\mathbb{E}\int_t^T \big\langle z(s),f(s)\big\rangle_H \,ds \\
\ns\ds
 =  \mathbb{E}\int_t^T
\big\langle \psi_1(s),y^t(s)\big\rangle_H
\,ds + \mathbb{E} \int_t^T\big\langle
\psi_2(s),Y^t(s)\big\rangle_{H^d} \,ds
+\mathbb{E}
\big\langle\eta,\xi^t\big\rangle_H.
\end{array}
\end{equation}
It is clear that $\xi^T=y_T$. Furthermore,
\begin{equation}\label{vw-th1-eq3}
\begin{array}{ll}\ds
|(y^t(\cdot), Y^t(\cdot),\xi^t)|_{
L^p_{\dbF}(\O; L^\infty(t,T;H)) \times L^2_{\dbF}(t,T;L^p(\O;H^d))\times L^p_{\cF_t}(\O;H)} \\
\ns\ds \leq C\left[
 |f(\cdot)|_{ L^1_{\dbF}(t,T;L^p(\O;H))}+|y_T|_{
L^p_{\cF_T}(\O;H)}\right],
\qq\q\forall\;t\in [\tau,T].
\end{array}
\end{equation}

\ms

{\bf Step 2.} Note that the  $(y^t(\cdot),
Y^t(\cdot))$ obtained in Step 1 may depend
on $t$. Now we  show the time consistency
of $(y^t(\cdot), Y^t(\cdot))$, that is, for
any $t_1$ and $t_2$ satisfying $\tau\leq
t_2 \leq t_1 \leq T$, it holds that
\begin{equation}\label{vw-th1-eq4}
\big(y^{t_2} (s,\o),Y^{t_2}
(s,\o)\big)=\big( y^{t_1}(s,\o),
Y^{t_1}(s,\o)\big),\qq \ae (s,\o) \in
[t_1,T]\times\O,
\end{equation}
by suitable choice of the $\eta$, $\psi_1$
and $\psi_2$ in \eqref{fsystem1}. In fact,
for any fixed $\varrho (\cdot)\in
L^q_{\dbF}(\O;L^1(t_1,T;H))$ and $\varsigma
(\cdot)\in L^2_{\dbF}(t_1,T;L^q(\O;H^d))$,
we choose first $t=t_1$, $\eta = 0$,
$\psi_1(\cdot)=\varrho (\cdot)$ and
$\psi_2(\cdot) = \varsigma (\cdot)$ in
\eqref{fsystem1}. From \eqref{vw-th1-eq2},
we obtain that
\begin{equation}\label{vw-th1-eq5}
 \ba{ll}\ds
\mathbb{E}\big\langle
z^{t_1}(T),y_T\big\rangle_H -
\mathbb{E}\int_{t_1}^T \big\langle
z^{t_1}(s),f(s)\big\rangle_H ds\\\ns\ds
 =  \mathbb{E}\int_{t_1}^T
\big\langle\varrho
(s),y^{t_1}(s)\big\rangle_H
ds+\mathbb{E}\int_{t_1}^T
\big\langle\varsigma
(s),Y^{t_1}(s)\big\rangle_{H^d} ds.
 \ea
\end{equation}
Then, we choose $t=t_2$, $\eta = 0$,
$\psi_1(t,\omega) = \chi_{[t_1,T]}(t)
\varrho (t,\omega)$ and $\psi_2(t,\omega) =
\chi_{[t_1,T]}(t) \varsigma (t,\omega)$ in
\eqref{fsystem1}. It follows from
\eqref{vw-th1-eq2} that
\begin{equation}\label{vw-th1-eq6}
\ba{ll}\ds \mathbb{E}\big\langle
z^{t_1}(T),y_T\big\rangle_H -
\mathbb{E}\int_{t_1}^T \big\langle
z^{t_1}(s),f(s)\big\rangle_H ds\\\ns\ds
 =  \mathbb{E}\int_{t_1}^T
\big\langle\varrho
(s),y^{t_2}(s)\big\rangle_H
ds+\mathbb{E}\int_{t_1}^T
\big\langle\varsigma
(s),Y^{t_2}(s)\big\rangle_{H^d} ds. \ea
\end{equation}
Combining  \eqref{vw-th1-eq5} and
\eqref{vw-th1-eq6}, we obtain that

$$
\begin{array}{ll}\ds
\mathbb{E}\int_{t_1}^T \big\langle\varrho
(s),y^{t_1}(s)-y^{t_2}(s)\big\rangle_H ds
+\mathbb{E}\int_{t_1}^T
\big\langle\varsigma
(s),Y^{t_1}(s)-Y^{t_2}(s)\big\rangle_{H^d}\,
ds=0,\\
\ns\ds\qq \qq\qq\forall\; \varrho
(\cdot)\in L^q_{\dbF}(\O;L^1(t_1,T;H)),\q
\varsigma (\cdot)\in
L^2_{\dbF}(t_1,T;L^q(\O;H^d)).
\end{array}
$$
This yields the desired equality
\eqref{vw-th1-eq4}.

Put
\begin{equation}\label{vw-th1-eq7}
 y(t,\o)=y^\tau(t,\o),\qq Y (t,\o)= Y^\tau(t,\o),\qq \forall\;(t,\o) \in
[\tau,T]\times\O.
\end{equation}
From  \eqref{vw-th1-eq4}, we see that
\begin{equation}\label{vw-th1-eq8}
 \big(y^t (s,\o),Y^t
(s,\o)\big)=\big( y(s,\o),Y (s,\o)\big),
\qq \ae(s,\o) \in [t,T]\times\O.
\end{equation}
Combining \eqref{vw-th1-eq2} and
\eqref{vw-th1-eq8}, we get that
\begin{equation}\label{vw-th1-eq9}
\begin{array}{ll}
\ds \mathbb{E}\big\langle
z(T),y_T\big\rangle_H - \mathbb{E}
\big\langle\eta,\xi^t \big\rangle_H\\
\ns\ds =\mathbb{E}\int_t^T \big\langle
z(s),f(s)\big\rangle_H ds+
\mathbb{E}\int_t^T \big\langle
\psi_1(s),y(s)\big\rangle_H ds +\mathbb{E}
\int_t^T \big\langle
\psi_2(s),Y(s)\big\rangle_{H^d} ds,\\\ns\ds
\qq\q\forall\; \big(\psi_1(\cdot),
\psi_2(\cdot),\eta\big)\in
L^q_{\dbF}(\O;L^1(t,T;H))\times
L^2_{\dbF}(t,T;L^q(\O;H^d))\times
L^q_{\cF_t}(\O;H).
\end{array}
\end{equation}

\ms

{\bf Step 3.} We show in this step that
$\xi^t$ has a c\`adl\`ag modification, and
$y(t,\o)=\xi^t(\o)$ for a.e. $(t,\o)\in
[\tau,T]\times\O$. The detail is lengthy
and very similar to Steps 3--4 in the proof
of \cite[Theorem 3.1]{LZ1}, and hence we
omit it here. This completes the proof of
Theorem \ref{vw-th1}.\endpf


\section{Well-posedness of the operator-valued BSEEs}
\label{sec well2}


In this section, we consider the
well-posedness of \eqref{bsystem2}.

In order to define the transposition
solution of \eqref{bsystem2}, for any $t\in
[\tau,T]$, we introduce the following two
(forward) stochastic evolution equations:
\begin{equation}\label{op-fsystem2}
\left\{
\begin{array}{ll}
\ds dx_1 = (A+J)x_1ds + u_1ds + Kx_1 dw(s) + v_1 dw(s) &\mbox{ in } (t,T],\\
\ns\ds x_1(t)=\xi_1
\end{array}
\right.
\end{equation}
and
\begin{equation}\label{op-fsystem3}
\left\{
\begin{array}{ll}
\ds dx_2 = (A+J)x_2ds + u_2ds + Kx_2 dw(s) + v_2 dw(s) &\mbox{ in } (t,T],\\
\ns\ds x_2(t)=\xi_2,
\end{array}
\right.
\end{equation}
where $\xi_i\in L^{2q}_{\cF_t}(\O;H)$,
$u_i\in L^{2q}_\dbF(\O;L^2(t,T;H))$,
$v_i\in L^2_\dbF(t,T;L^{2q}(\O;H^d))$ and
$i=1,2$. Also, we need to introduce the
solution space for the equation
\eqref{bsystem2}. Put
 $$
\begin{array}{ll}\ds
L^{p}_{\dbF,w}(\O; D([\tau,T];\cL(H))\\
\ns\ds\= \Big\{P(\cd,\cd)\;\Big|\;
P(\cd,\cd)\in
\cL_{pd}\big(L^{2q}_{\dbF}(\O;L^{2}(\tau,T;H)),\;L^{\frac{2p}{p+1}}_{\dbF}(\O;L^2(\tau,T;H))\big), \mbox{ and for every }\\
\ns\ds\qq\qq\q  t\in[\tau,T]\hb{ and
}\xi\in
L^{2q}_{\cF_t}(\O;H),\,P(\cd,\cd)\xi\in
L^{\frac{2p}{p+1}}_{\dbF}(\O;D([t,T];H))\\
\ns\ds\qq\qq\q  \mbox{and }
|P(\cd,\cd)\xi|_{L^{\frac{2p}{p+1}}_{\dbF}(\O;D([t,T];H))}
\leq C|\xi|_{L^{2q}_{\cF_t}(\O;H)} \Big\}
\end{array}
 $$
and
$$
L^2_{\dbF,w}(\tau,T;L^{p}(\O;\cL(H)))^d\=
\big[\cL_{pd}\big(L^{r}_{\dbF}(\tau,T;L^{2q}(\O;H^d)),
\;L^{\frac{2r}{2+r}}_{\dbF}(\tau,T;
L^{\frac{2p}{p+1}}(\O;H))\big)\big],
$$
where $r\geq 2$. The transposition solution
to the equation \eqref{bsystem2} is defined
as follows:

\begin{definition}\label{OP-def1}
We call $(P(\cd),Q(\cd))\in
L^{p}_{\dbF,w}(\O;D([\tau,T];
\cL(H)))\times
L^2_{\dbF,w}(\tau,T;L^{p}(\O;\cL(H)^d))$ a
transposition solution to  the equation
\eqref{bsystem2} if for any $t\in
[\tau,T]$, $\xi_1,\xi_2\in
L^{2q}_{\cF_t}(\O;H)$, $u_1(\cd),
u_2(\cd)\in L^{2q}_{\dbF}(\O;L^2(t,T;H))$
and $v_1(\cd),v_2(\cd)\in
L^2_{\dbF}(t,T;L^{2q}(\O; H^d))$, it holds
that
 \begin{equation}\label{OP-def1-eq1}
\begin{array}{ll}\ds
\dbE \big\langle P_T
x_1(T),x_2(T)\big\rangle_{H}
 - \dbE\int_t^T \big\langle F(s)x_1(s),x_2(s)\big\rangle_{H}ds\\
\ns\ds = \dbE \big\langle
P(t)\xi_1,\xi_2\big\rangle_{H} +
\dbE\int_t^T \big\langle
P(s)u_1(s),x_2(s)\big\rangle_{H} ds +
\dbE\int_t^T \big\langle P(s)x_1(s),u_2(s)\big\rangle_{H} ds \\
\ns\ds \q + \dbE\int_t^T \big\langle P(s)
K(s)x_1(s), v_2(s)\big\rangle_{H^d} ds +
\dbE\int_t^T \big\langle P(s)v_1(s),
K(s)x_2(s)+v_2(s)\big\rangle_{H^d} ds\\
\ns\ds \q  +  \dbE\int_t^T \big\langle
Q(s)v_1(s),x_2(s)\big\rangle_{H} ds +
\dbE\int_t^T \big\langle
x_1(s),Q(s)^*v_2(s)\big\rangle_{H} ds,
\end{array}
\end{equation}
where, $x_1(\cd)$ and $x_2(\cd)$ solve
\eqref{op-fsystem2} and
\eqref{op-fsystem3}, respectively.
\end{definition}

The well-posedness of the equation
\eqref{bsystem2} in the sense of Definition
\ref{OP-def1} is still open. However, we
can show the following uniqueness result
for the transposition solution to
\eqref{bsystem2}.

\begin{theorem}\label{OP-U-th}
Assume that $J, K, F$ and $P_T$ satisfy
\eqref{JKFP}. Then the equation
\eqref{bsystem2} admits at most one
transposition solution $(P(\cd),Q(\cd))\in
L^{p}_{\dbF,w}(\O;D([\tau,T];
\cL(H)))\times
L^2_{\dbF,w}(\tau,T;L^{p}(\O;\cL(H)^d))$.
\end{theorem}

{\it Proof}\,: Assume that $(\overline
P(\cd),\overline Q(\cd))\in
L^{p}_{\dbF,w}(\O;D([\tau,T];
\cL(H)))\times
L^2_{\dbF,w}(\tau,T;L^{p}(\O;\cL(H)^d))$ is
another transposition solution to
\eqref{bsystem2}. Then, by Definition
\ref{OP-def1}, it follows that, for any
$t\in [\tau,T]$,
\begin{equation}\label{OP-U-th-eq1}
\begin{array}{ll}
\ds 0=\mE\big\langle \big[\,\overline P(t)
- P(t)\big] \xi_1,\xi_2 \big\rangle_{H} +
\mE \int_t^T \big\langle \big[\,\overline
P(s) - P(s)\big]u_1(s),
x_2(s)\big\rangle_{H}ds \\
\ns\ds \qq + \mE \!\int_t^T\! \big\langle
\big[\,\overline P(s)\! -
\!P(s)\big]x_1(s), u_2(s)\big\rangle_{H}ds
\! +\! \mE\! \int_t^T\!\! \big\langle
\big[\,\overline P(s)\! -\! P(s)\big] K
(s)x_1 (s), v_2 (s)\big\rangle_{H^d}ds \\
\ns\ds \qq  +  \mE \int_t^T \big\langle
\big[\,\overline P(s) - P(s)\big]v_1
(s), K (s)x_2 (s)+v_2(s)  \big\rangle_{H^d}ds \\
\ns\ds \qq + \mE \!\int_t^T\!
\big\langle\big[\,\overline Q(s)\! -\!
Q(s)\big] v_1(s), x_2(s)\big\rangle_{H}ds\!
+\! \mE\! \int_t^T\!\! \big\langle x_1(s),
\big[\,\overline Q(s)^*\! -\!
Q(s)^*\big]v_2(s) \big\rangle_{H}ds.
\end{array}
\end{equation}
Choosing $u_1=u_2=0$ and $v_1=v_2=0$ in the
equation \eqref{op-fsystem2} and the
equation \eqref{op-fsystem3}, respectively,
by \eqref{OP-U-th-eq1}, we obtain that
$$
0=\mE\big\langle \big[\,\overline P(t) -
P(t)\big] \xi_1,\xi_2
\big\rangle_{H},\qq\forall\;t\in [\tau,T],
\forall\; \xi_1,\; \xi_2\in
L^{2q}_{\cF_t}(\O;H).
$$
Hence, we find that $\overline
P(\cd)=P(\cd)$. By this, it is easy to see
that for any $t\in [\tau,T]$,
\begin{equation}\label{OP-U-th-eq2}
0=\mE \int_t^T \big\langle\big[\,\overline
Q(s) - Q(s)\big] v_1(s),
x_2(s)\big\rangle_{H}ds  + \mE \int_t^T
\big\langle x_1(s), \big[\,\overline Q(s)^*
- Q(s)^*\big]v_2(s) \big\rangle_{H}ds.
\end{equation}
Choosing $t=\tau$, $\xi_2=0$ and $v_2=0$ in
\eqref{op-fsystem3}, we see that
\eqref{OP-U-th-eq2} becomes
\begin{equation}\label{OP-U-th-eq3}
0=\mE \int_\tau^T
\big\langle\big[\,\overline Q(s) -
Q(s)\big] v_1(s), x_2(s)\big\rangle_{H}ds.
\end{equation}

Similar to the proof of \cite[Theorem
4.1]{LZ1}, one can show that the set
$$
\Xi\=\left\{x_2(\cd)\;\left|\; x_2(\cd)\hb{
solves } \eqref{op-fsystem3} \hb{ with
}t=\tau,\;\xi_2=0,\; v_2=0\hb{ and some
}u_2\in L^{2q}_{\dbF}(\Omega;L^2(\tau,T;H))
\right.\right\}
$$
is dense in
$L^4_{\dbF}(0,T;L^{2q}(\Omega;H))$. By this
fact, and noting \eqref{OP-U-th-eq3}, we
see that
$$
\big[\,\overline Q(\cd) - Q(\cd)\big]
v_1(\cd)=0,\qq\forall\; v_1(\cdot)\in
L^4_{\dbF}(0,T;L^{2q}(\Omega;H^d)).
$$
Therefore, we find that we see that for all
$r\geq 2$,
$$
\big[\,\overline Q(\cd) - Q(\cd)\big]
v_1(\cd)=0,\qq\forall\; v_1(\cdot)\in
L^r_{\dbF}(0,T;L^{2q}(\Omega;H^d)).
$$
Hence $\overline Q(\cd) = Q(\cd)$. This
completes the proof of Theorem
\ref{OP-U-th}.
\endpf

\bs

Further, we have the following
well-posedness result for \eqref{bsystem2}
in a special case.

\begin{theorem}\label{Op-th1}
If $H$ is a separable Hilbert space,
$P_T\in L^p_{\cF_T}(\O;\cL_2(H))$, $F\in
L^1_\dbF(\tau,T;L^p(\O;\cL_2(H)))$, $J \in
L^4_\dbF(\tau,T; L^\infty(\O; \cL(H)))$ and
$K\in L^4_\dbF(\tau,T;
L^\infty(\O;\cL(H)^d))$, then
\eqref{bsystem2} admits a unique
transposition solution $ \big(P(\cd),
Q(\cd)\big) \in
L^{p}_\dbF(\O;D([\tau,T];\cL_2(H)))\times
L^2_{\dbF}(\tau,T;L^p(\O;\cL_2(H)^d))$.
Furthermore,
\begin{equation}\label{OP-th1-eq1}
|(P, Q)|_{
L^{p}_\dbF(\O;D([\tau,T];\cL_2(H)))\times
L^2_\dbF(\tau,T;L^{p}(\O;\cL_2(H)^d))}\leq
C\big(
|F|_{L^1_\dbF(\tau,T;L^p(\O;\cL_2(H)))} +
|P_T|_{L^p_{\cF_T}(\O;\cL_2(H))}\big).
\end{equation}
\end{theorem}

{\it Proof}\,: The proof is very similar to
that for \cite[Theorem 4.2]{LZ1}, and hence
we only give below a sketch.

First, we define a family of operators
$\{\cT (t)\}_{t\geq 0}$ on $\cL_2(H)$ as
follows:
$$
\cT (t)O = S(t)OS^*(t), \q\forall\; O\in
\cL_2(H).
$$
Then, $\{\cT (t)\}_{t\geq 0}$ is a
$C_0$-semigroup on $\cL_2(H)$. Denote by
$\cA$ the infinitesimal generater of $\{\cT
(t)\}_{t\geq 0}$. We consider the following
$\cL_2(H)$-valued BSEE:
\begin{equation}\label{OP-HS-equ1}
\left\{
\begin{array}{ll}\ds
dP = -\cA^* Pdt + f(t,P,Q)dt + Qdw &\mbox{
in
}[\tau,T),\\
\ns\ds P(T)=P_T,
\end{array}
\right.
\end{equation}
where
 $
f(t,P,Q) = -J^*P-PJ-K^*PK - K^*Q - QK + F$.
Since $J\in L^4_\dbF(\tau,T;
L^\infty(\O;\cL(H)))$, $K\in
L^4_\dbF(\tau,T; L^\infty(\O;\cL(H)^d))$
and $F\in
L^1_\dbF(\tau,T;L^p(\O;\cL_2(H)))$, we see
that $f(\cd,\cd,\cd)$ satisfies
\eqref{Lm1}. Since $\cL_2(H)$ is a Hilbert
space, by Theorem \ref{vw-th1}, one can
find a pair $(P, Q) \in
L^{p}_\dbF(\O;D([\tau,T];\cL_2(H)))\times
L^2_\dbF(\tau,T;L^{p}(\O;\cL_2(H)^d))$
solving the equation \eqref{OP-HS-equ1} in
the sense of Definition \ref{definition1}.
Further, $(P,Q)$ satisfies
\eqref{OP-th1-eq1}.

Next, denote by $O(\cd)$ the tensor product
of $x_1(\cd)$ and $x_2(\cd)$, solutions to
\eqref{op-fsystem2} and
\eqref{op-fsystem3}, respectively. As
usual, $O(s,\o)x = \langle x,x_1(s,\o)
\rangle_H x_2(s,\o)$ for a.e. $(s,\o)\in
[t,T]\times \O$ and $x\in H$. Hence,
$O(s,\o)\in \cL_2(H)$  for a.e. $(s,\o)\in
[t,T]\times \O$. It can be proved that
\begin{equation}\label{OP-U-eq3}
\left\{\3n
\begin{array}{ll}
\ds dO(s)  = \cA O(s) ds + uds + vdw(s) &\mbox{ in } (t,T],\\
\ns\ds O(t) = \xi_1\otimes \xi_2,
\end{array}
\right.
\end{equation}
where
 $$
\left\{
\begin{array}{ll}\ds
u  = J O (\cd) + O(\cd)  J^* + u_1\otimes
x_2  + x_1 \otimes u_2 + K O(\cd)  K^*  +
(K x_1 ) \otimes v_2 + v_1 \otimes (K x_2)
+ v_1 \otimes
v_2 , \\
\ns\ds v  = K O (\cd) + O (\cd) K^* + v_1
\otimes x_2  + x_1  \otimes v_2.
\end{array}
\right.
 $$

Noting that $\big(P(\cd),Q(\cd)\big)$
solves the equation \eqref{OP-HS-equ1} in
the transposition sense and by
\eqref{OP-U-eq3}, we have that
\begin{equation}\label{OP-U-eq4}
\begin{array}{ll}\ds
\dbE \big\langle
O(T),P_T\big\rangle_{\cL_2(H)}
- \dbE\int_t^T \big\langle O(s),f(s,P(s),Q(s) )\big\rangle_{\cL_2(H)}ds\\
\ns\ds = \dbE
\big\langle\xi_1\otimes\xi_2,P(t)\big\rangle_{\cL_2(H)}
+ \dbE\int_t^T \big\langle
u(s),P(s)\big\rangle_{\cL_2(H)} ds +
\dbE\int_t^T \big\langle
v(s),Q(s)\big\rangle_{\cL_2(H)^d} ds.
\end{array}
\end{equation}

Finally, by (\ref{OP-U-eq4}) and some
direct computation, one can show that
$\big(P(\cd),Q(\cd)\big)$ satisfies
\eqref{OP-def1-eq1}, and therefore it is a
transposition solution to the equation
\eqref{bsystem2} (in the sense of
Definition \ref{OP-def1}). The uniqueness
of $\big(P(\cd),Q(\cd)\big)$ follows from
Theorem \ref{OP-U-th}.
\endpf

\ms

Since we are not able to prove the
well-posedness of the equation
\eqref{bsystem2} in the sense of Definition
\ref{OP-def1} at this moment, we need to
introduce a weaker notion of solution,
i.e., relaxed transposition solution to
this equation. For this purpose, we write
 $$
\3n\begin{array}{ll}\ds
\cQ^p[\tau,T]\\
\ns\ds\=\Big\{\big(Q^{(\cd)},\widehat
Q^{(\cd)}\big)\;\Big|\;
Q^{(\cd)}=(Q^{1,(\cd)},\cds,Q^{d,(\cd)}),\,\widehat
Q^{(\cd)}=(\widehat
Q^{1,(\cd)},\cds,\widehat Q^{d,(\cd)}).
\mbox{ For arbitrary } t\in [\tau,T]
 \\
\ns\ds \q  \mbox{and }\, i=1,\cds,d,
\,\mbox{ both }Q^{i,(t)}\mbox{ and
}\widehat Q^{i,(t)}\mbox{
are bounded linear operators from } L^{2q}_{\cF_t}(\O;H)\\
\ns\ds \q \times
L^{2q}_\dbF\!(\O;L^2(t,T;\!H))\!\times\!
L^2_\dbF(t,T;L^{2q}(\O;\!H)) \mbox{ to }
L^{\frac{2p}{p+1}}_\dbF\!(\O;L^{2}(t,T;\!H))\!
\mbox{ and
}\!Q^{(t)}(0,0,\cd)^*\!=\!\widehat
Q^{(t)}(0,0,\cd)\!\Big\}.
\end{array}
 $$
For $\big(Q^{(\cd)},\widehat
Q^{(\cd)}\big)\in \cQ^p[\tau,T]$, put
$$
\begin{array}{ll}\ds
\big|\big(Q^{(\cd)},\widehat
Q^{(\cd)}\big)\big|_{\cQ^p[\tau,T]}\\
\ns\ds\=\sum_{i=1}^d\sup_{t\in
[\tau,T]}\big|\!\big|\big(Q^{i,(t)},\widehat
Q^{i,(t)}\big)\big|\!\big|_{\big(\cL(L^{2q}_{\cF_t}(\O;H)\times
L^{2q}_\dbF(\O;L^2(t,T;H))\times
L^2_\dbF(t,T;L^{2q}(\O;H)),\;
L^{\frac{2p}{p+1}}_\dbF(\O;L^{2}(t,T;H))\big)^2}.
\end{array}
$$

The relaxed transposition solution to
\eqref{bsystem2} is defined as follows:

\begin{definition}\label{OP-def3}
We call $\big(P(\cd),Q^{(\cd)},\widehat
Q^{(\cd)}\big)\in
L^{p}_{\dbF,w}(\O;D([\tau,T];
\cL(H)))\times \cQ^p[\tau,T]$ a relaxed
transposition solution to the equation
\eqref{bsystem2} if for any $t\in
[\tau,T]$, $\xi_1,\xi_2\in
L^{2q}_{\cF_t}(\O;H)$, $u_1(\cd),
u_2(\cd)\in L^{2q}_{\dbF}(\O;L^2(t,T;H))$
and $v_1(\cd),v_2(\cd)\in
L^2_{\dbF}(t,T;L^{2q}(\O; H^d))$, it holds
that
\begin{equation}\label{6.18eq1}
\begin{array}{ll}
\ds \mE\big\langle P_T x_1(T), x_2(T)
\big\rangle_{H} - \mE \int_t^T \big\langle
F(s) x_1(s), x_2(s) \big\rangle_{H}ds\\
\ns\ds =\mE\big\langle P(t) \xi_1,\xi_2
\big\rangle_{H} + \mE \int_t^T \big\langle
P(s)u_1(s), x_2(s)\big\rangle_{H}ds + \mE
\int_t^T \big\langle P(s)x_1(s),
u_2(s)\big\rangle_{H}ds \\
\ns\ds \q  + \mE \int_t^T \big\langle P(s)K
(s)x_1 (s), v_2 (s)\big\rangle_{H^d}ds +
\mE
\int_t^T \big\langle  P(s)v_1 (s), K (s)x_2 (s)+ v_2(s)\big\rangle_{H^d}ds\\
\ns\ds \q + \mE \int_t^T \big\langle
v_1(s), \widehat
Q^{(t)}(\xi_2,u_2,v_2)(s)\big\rangle_{H^d}ds+
\mE \int_t^T \big\langle
Q^{(t)}(\xi_1,u_1,v_1)(s), v_2(s)
\big\rangle_{H^d}ds,
\end{array}
\end{equation}
where, $x_1(\cd)$ and $x_2(\cd)$ solve the
equations \eqref{op-fsystem2} and
\eqref{op-fsystem3}, respectively.
\end{definition}

We have the following well-posedness result
for  the equation \eqref{bsystem2}.

\begin{theorem}\label{OP-th2}
Assume that $H$ is a separable Hilbert
space, and $L^p_{\cF_T}(\O;\dbC)$ ($1\leq p
< \infty$) is a separable Banach space.
Then, for any $J$, $K$, $F$ and $P_T$
satisfying (\ref{JKFP}), the equation
\eqref{bsystem2} admits one and only one
relaxed transposition solution
$\big(P(\cd),Q^{(\cd)},\widehat
Q^{(\cd)}\big)$. Furthermore,
 $$
\begin{array}{ll}\ds
\q |P|_{ L^{p}_{\dbF,w}(\O;D([\tau,T];
\cL(H)))} + \big|\big(Q^{(\cd)},\widehat
Q^{(\cd)}\big)\big|_{\cQ^p[\tau,T]}\leq
C\big[
|F|_{L^1_\dbF(\tau,T;\;L^p(\O;\cL(H)))} +
|P_T|_{L^p_{\cF_T}(\O;\;\cL(H))}\big].
\end{array}
 $$
\end{theorem}

{\it Proof}\,: The proof of this theorem is
very lengthy and technical, and it is very
similar to that of \cite[Theorem 6.1]{LZ1}.
Hence, we only give here a sketch.

In Theorem \ref{OP-U-th}, we have obtained
the well-posedness of \eqref{bsystem2} with
$P_T\in L^p_{\cF_T}(\O;\cL_2(H))$ and $F\in
L^1_\dbF(\tau,T;L^p(\O;\cL_2(H)))$. Noting
that $\cL_2(H)$ is dense in the space
$\cL(H)$ (with the usual strong operator
topology), we may approximate the general
data $P_T\in L^p_{\cF_T}(\O;\cL(H))$ and
$F\in L^1_\dbF(\tau,T;L^p(\O;\cL(H)))$ by
$\{P_T^m\}_{m=1}^\infty\subset
L^p_{\cF_T}(\O;\cL_2(H))$ and
$\{F^m\}_{m=1}^\infty\!\subset\!
L^1_\dbF(\tau,T;L^p(\O;\cL_2(H)))$,
respectively. Denote by
$(P^m(\cd),Q^m(\cd))$ the corresponding
solution to \eqref{bsystem2} with $P_T$ and
$F$ replaced respectively by $P_T^m$ and
$F^m$. Then, we obtain the desired $P(\cd)$
as the weak limit of
$\{P^n(\cd)\}_{m=1}^\infty$, and
$(Q^{(\cd)},\widehat Q^{(\cd)})$ as the
weak limit of
$\{(Q^m(\cd),Q^m(\cd)^*)\}_{m=1}^\infty$.
The most difficult part is to show that
$\{P^n(\cd)\}_{m=1}^\infty$ and
$\{(Q^m(\cd),Q^m(\cd)^*)\}_{m=1}^\infty$
converge respectively to some elements in
$L^{p}_{\dbF,w}(\O;D([\tau,T];
\cL(H)))\times \cQ^p[\tau,T]$, in some weak
sense. All of these are guaranteed by some
Banach-Alaoglu-type theorems established in
\cite{LZ1}.
\endpf


\section{A regularity property for relaxed transposition solutions to the operator-valued BSEEs}\label{as4}
\label{sec regular}

In this section, we shall derive a
regularity property for relaxed
transposition solutions to the equation
\eqref{bsystem2}. This property will play
key roles in the proof of our general
Pontryagin-type stochastic maximum
principle, presented in Section \ref{max}.
To simplify the notations, we assume that
$d=1$ in this section.

We need some preliminaries. First of all,
as an immediate consequence of the
well-posedness result for
\eqref{op-fsystem3}, it is easy to prove
the following result.

\begin{lemma}\label{4.20-lm1}
If $u_2=v_2=0$ in the equation
\eqref{op-fsystem3}, then for each
$t\in[0,T]$, there exists an operator
$U(\cd,t)\in \cL\big(L^p_{\cF_t}(\O;H),
C_\dbF([t,T];L^p(\O;H))\big)$ such that the
corresponding solution to this equation can
be represented as $x_2(\cd) =
U(\cd,t)\xi_2$. Further, for any $t\in
[0,T]$, $\e>0$ and $\xi\in
L^p_{\cF_t}(\O;H)$, there is a $\d>0$ such
that for all $s\in [t,t+\d]$,
$$
|U(s,t)\xi-\xi|_{L^p_{\cF_T}(\O;H)}<\e.
$$
\end{lemma}

Next, we recall the following known result.

\begin{lemma} {\rm (\cite[Corollary 5.1]{LZ1})} \label{cor1}
Let $X$ and $Y$ be respectively a separable
and a reflexive Banach space, and let
$L^q_{\cF_T}(\O)$, with $1\leq q <\infty$,
be separable. Let $1<q_1,q_2 < \infty$.
Assume that $\{\cG_n\}_{n=1}^\infty$ is a
sequence of uniformly bounded, pointwisely
defined linear operators from $X$ to
$L^{q_1}_{\dbF}(0,T;L^{q_2}(\O;Y))$. Then,
there exist a subsequence
$\{\cG_{n_k}\}_{k=1}^\infty\subset
\{\cG_n\}_{n=1}^\infty$ and an
$\cG\in\cL_{pd}\big(X,\;
L^{q_1}_{\dbF}(0,T;L^{q_2}(\O;Y))\big)$
such that
$$
\cG x = {\mbox{\rm(w)-}}\lim_{k\to\infty}
\cG_{n_k}x \ \mbox{ in }
L^{q_1}_{\dbF}(0,T;L^{q_2}(\O;Y)), \qq
\forall\; x\in X.
$$
Moreover,
 $\ds|\!|\cG|\!|_{\cL(X,
L^{q_1}_{\dbF}(0,T;L^{q_2}(\O;Y)))}\le
\sup_{n\in\dbN}|\!|\cG_n|\!|_{\cL(X, \
L^{q_1}_{\dbF}(0,T;L^{q_2}(\O;Y))}$.
\end{lemma}

Further, let $\{\D_n\}_{n=1}^\infty$ be a
sequence of partitions of $[0,T]$, that is,

$$
\D_n\= \Big\{t_i^n\;\Big|\;i=0,1,\cdots,n,
\hb{ and }0=t_0^n < t_1^n < \cds < t_{n}^n
=T\Big\}
$$
such that $\D_n\subset \D_{n+1}$ and
$\d(\D_n)\= \max_{0\leq i\leq n-1}
(t_{i+1}^n - t_{i}^n)\to 0$ as
$n\to\infty$. We introduce the following
subspaces of $L^2_\dbF(0,T;L^{2q}(\O;H))$:
\begin{equation}\label{cH}
\cH_n=\Big\{\sum_{i=0}^{n-1}
\chi_{[t_i^n,t_{i+1}^n)}(\cd)U(\cd,t_i^n)h_i\;\Big|\;
h_i\in L^{2q}_{\cF_{t_i^n}}(\O;H)\Big\}.
\end{equation}
Here $U(\cd,\cd)$ is the operator
introduced in Lemma \ref{4.20-lm1}. We have
the following result.
\begin{proposition}\label{5.13-prop1}
The set $\bigcup_{n=1}^\infty \cH_n$ is
dense in $L^2_\dbF(0,T;L^{2q}(\O;H))$.
\end{proposition}
{\it Proof}\,: We introduce the following
subspace of $L^2_\dbF(0,T;L^{2q}(\O;H))$:
\begin{equation}\label{cH1}
\wt\cH_n=\Big\{\sum_{i=0}^{n-1}
\chi_{[t_i^n,t_{i+1}^n)}(\cd)h_i^n\;\Big|\;
h_i^n\in L^{2q}_{\cF_{t_i^n}}(\O;H)\Big\}.
\end{equation}
It is clear that
$\bigcup_{n=1}^\infty\wt\cH_n$ is dense in
$L^2_\dbF(0,T;L^{2q}(\O;H))$.

For any $n\in \dbN$ and $h_i^n\in
L^{2q}_{\cF_{t_i^n}}(\O;H)$,
$i\in\{0,1,\cdots,n-1\}$, write
 $\tilde v_n \=
\sum_{i=0}^{n-1}
\chi_{[t_i^n,t_{i+1}^n)}(\cd)h_i^n$.
Clearly, $\tilde v_n \in \wt H_n$. We claim
that for any $\e>0$, there exist an $m\in
\dbN$ and a $v_m\in \cH_m$ such that
\bel{zc1} \big|\tilde v_n -
v_m\big|_{L^2_\dbF(0,T;L^{2q}(\O;H))}<\e.
\ee

Indeed, by Lemma \ref{4.20-lm1}, for each
$h_i^n$, there is a $\d_i^n>0$ such that
for all $t\in [t_i^n,T-\d_i^n)$ and $s\in
[t,t+\d_i^n]$, it holds that
\begin{equation}\label{5.13-eq1}
\Big|U(s,t)h_i^n-h_i^n\Big|_{L^{2q}_{\cF_T}(\O;H)}<\frac{\e}{\sqrt{T}}.
\end{equation}
Now we choose a partition $\D_m$ of $[0,T]$
such that $\D_n\subset\D_m$ and
$\max_{0\leq j\leq m-1} \{t_{j+1}-t_j\}\leq
\min_{0\leq i\leq n-1}\{\d_i^n\}$. Let
$$
v_m=\sum_{j=0}^{m-1}\chi_{[t_j^m,t_{j+1}^m)}(\cd)U(\cd,t_j^m)h_j^m,
$$
where $h_j^m=h_i^n$ whenever
$[t_j^m,t_{j+1}^m)\subset
[t_i^n,t_{i+1}^n)$. From \eqref{5.13-eq1},
we find that
$$
\begin{array}{ll}\ds
\big|\tilde v_n -
v_m\big|_{L^2_\dbF(0,T;L^{2q}(\O;H))}=\Big|\sum_{j=0}^{m-1}
\chi_{[t_j^m,t_{j+1}^m)}(\cd)U(\cd,t_j^m)h_j^m
- \sum_{i=0}^{n-1}
\chi_{[t_i^n,t^n_{i+1})}(\cd)h_i^n\Big|_{L^2_\dbF(0,T;L^{2q}(\O;H))}\\
\ns\ds = \Big|\sum_{j=0}^{m-1}
\chi_{[t_j^m,t_{j+1}^m)}(\cd)\[U(\cd,t_j^m)h_j^m-h_j^m\]\Big|_{L^2_\dbF(0,T;L^{2q}(\O;H))}<
\sqrt{T}\frac{\e}{\sqrt{T}}=\e.
\end{array}
$$
This proves \eqref{zc1}. Hence,
$\bigcup_{n=1}^\infty \cH_n$ is dense in
$L^2_\dbF(0,T;L^{2q}(\O;H))$.
\endpf

Our regularity result for solutions to
\eqref{bsystem2} can be stated as follows.

\begin{theorem}\label{10.1th}
Suppose that the assumptions in Theorem
\ref{OP-th2} hold and let
$(P(\cd),Q^{(\cd)},\widehat Q^{(\cd)})$ be
the relaxed transposition solution to the
equation \eqref{bsystem2}. Then,  there
exist an $n\in\dbN$ and two pointwisely
defined linear operators $Q^n$ and
$\widehat Q^n$, both of which are from
$\cH_n$ to
$L^{2}_\dbF(0,T;L^{\frac{2p}{p+1}}(\O;H))$,
such that, for any $\xi_1,\xi_2\in
L^{2q}_{\cF_0}(\O;H)$, $u_1(\cd),
u_2(\cd)\in L^{2q}_\dbF(\O;L^2(0,T;H))$ and
$v_1(\cd),v_2(\cd)\in \cH_{n}$, it holds
that
\begin{equation}\label{10.9eq2}
\begin{array}{ll}
\ds \mE \int_{0}^T \big\langle v_1(s),
\widehat Q^{(0)}(\xi_2,u_2,v_2)(s)
\big\rangle_{H}ds + \mE \int_{0}^T
\big\langle
 Q^{(0)}(\xi_1,u_1,v_1) (s),
v_2(s) \big\rangle_{H}ds \\
\ns\ds =\mE \int_{0}^T
 \[\big\langle
\big(Q^n v_1\big)(s), x_2 (s)
\big\rangle_{H}+\big\langle x_1 (s),
\big(\widehat Q^n v_2\big)(s)
\big\rangle_{H}\]ds,
\end{array}
\end{equation}
where, $x_1(\cd)$ and $x_2(\cd)$ solve
accordingly \eqref{op-fsystem2} and
\eqref{op-fsystem3} with $t=0$. Further,
there is a positive constant $C$,
independent of $n$, such that
\begin{equation}\label{4.14-eq3}
\big|Q^n
v_1\big|_{L^{2}_\dbF(0,T;L^{\frac{2p}{p+1}}(\O;H))}
+ \big|\widehat Q^n
v_2\big|_{L^{2}_\dbF(0,T;L^{\frac{2p}{p+1}}(\O;H))}\leq
C\big(|\tilde
v_1|_{L^2_\dbF(0,T;L^{2q}(\O;H))}+|\tilde
v_2|_{L^2_\dbF(0,T;L^{2q}(\O;H))}\big),
\end{equation}
where
$$
\tilde v_1 = \sum_{i=0}^{n-1}
\chi_{[t^n_i,t^n_{i+1})}(\cd) h_i \q\mbox{
for }\; v_1 = \sum_{i=0}^{n-1}
\chi_{[t^n_i,t^n_{i+1})}(\cd)U(\cd,t_i) h_i
$$
and
$$
\tilde v_2 = \sum_{j=0}^{n-1}
\chi_{[t^n_j,t^n_{j+1})}(\cd) h_j  \q\mbox{
for }\; v_2 = \sum_{j=0}^{n-1}
\chi_{[t^n_j,t^n_{j+1})}(\cd)U(\cd,t_j)h_j.
$$
\end{theorem}

{\it Proof:}\, Let $\{e_m\}_{m=1}^\infty$
be an orthonormal basis of $H$ and
$\{\G_m\}_{m=1}^\infty$ be the standard
projection operator from $H$ onto its
subspace $\span \{e_1,e_2,\cdots,e_m\}$.
Let $P_T^m=\G_mP_T\G_m$ and $F_m(\cd)=\G_m
F(\cd)\G_m$. Clearly, $P_T^m\in
L^2_{\cF_T}(\O;\cL_2(H))$ and $F_m\in
L^1_\dbF(0,T;L^2(\O;\cL_2(H)))$. By Theorem
\ref{Op-th1}, the equation \eqref{bsystem2}
with $P_T$ and $F$ replaced respectively by
$P_T^m$ $F_m$ admits a unique transposition
solution $(P^m,Q^m)\in
L^{p}_\dbF(\O;D([\tau,T];\cL_2(H)))\times
L^2_{\dbF}(\tau,T;L^p(\O;\cL_2(H)))$ such
that
\begin{equation}\label{4.14-eq2}
\begin{array}{ll}\ds
\dbE \big\langle P_T^m
x_1(T),x_2(T)\big\rangle_{H}
- \dbE\int_t^T \big\langle F_m(s)x_1(s),x_2(s)\big\rangle_{H}ds\\
\ns\ds = \dbE \big\langle
P^m(t)\xi_1,\xi_2\big\rangle_{H} +
\dbE\int_t^T \big\langle
P^m(s)u_1(s),x_2(s)\big\rangle_{H} ds +
\dbE\int_t^T \big\langle P^m(s)x_1(s),u_2(s)\big\rangle_{H} ds \\
\ns\ds \q + \dbE\int_t^T \big\langle P^m(s)
K(s)x_1(s), v_2(s)\big\rangle_{H} ds +
\dbE\int_t^T \big\langle P^m(s)v_1(s),
K(s)x_2(s)+v_2(s)\big\rangle_{H} ds\\
\ns\ds \q  + \dbE\int_t^T \big\langle
Q^m(s)v_1(s),x_2(s)\big\rangle_{H} ds +
\dbE\int_t^T \big\langle
Q^m(s)x_1(s),v_2(s)\big\rangle_{H} ds.
\end{array}
\end{equation}
Here, $x_1(\cd)$ and $x_2(\cd)$ solve
\eqref{op-fsystem2} and
\eqref{op-fsystem3}, respectively.

For any $i\in \{1,2,\cdots,n-1\}$,
$\xi_1\in L^{2q}_{\cF_{t^n_i}}(\O;H)$ and
$v_2\in L^2_\dbF(t_i^n,T;L^{2q}(\O;H))$,
letting $u_1 = 0$ and $v_1=0$ in the
equation \eqref{op-fsystem2}, and letting
$\xi_2=0$ and $u_2=0$ in the equation
\eqref{op-fsystem3}, by \eqref{4.14-eq2}
with $t=t_i^n$, we find that
\begin{equation}\label{6.16eq1}
\begin{array}{ll}
\ds \mE\big\langle P^m_T x_1(T), x_2(T)
\big\rangle_{H} - \mE \int_{t_i^n}^T
\big\langle
F_m(s) x_1(s), x_2(s) \big\rangle_{H}ds\\
\ns\ds = \mE \int_{t_i^n}^T \big\langle
P^m(s)K(s)x_1(s),
v_2(s)\big\rangle_{H}ds+\mE \int_{t_i^n}^T
\big\langle Q^m(s)U(s,t^n_i)\xi_1, v_2(s)
\big\rangle_{H}ds.
\end{array}
\end{equation}
For these data $\xi_1$, $u_1$, $v_1$,
$\xi_2$, $u_2$ and $v_2$, from the
variational equality \eqref{4.14-eq2} with
$t=t_{i+1}^n$, we obtain that
\begin{equation}\label{6.16eq1xx}
\begin{array}{ll}
\ds \mE\big\langle P^{m}_T x_1(T), x_2(T)
\big\rangle_{H} - \mE \int_{t_{i+1}^n}^T
\big\langle
F_m(s) x_1(s), x_2(s) \big\rangle_{H}ds\\
\ns\ds = \mE \big\langle
P^{m}(t_{i+1}^n)x_1(t_{i+1}^n),
x_2(t_{i+1}^n)\big\rangle_{H}  + \mE
\!\int_{t_{i+1}^n}^T \! \!\big\langle
P^{m}(s)K(s)x_1(s),
v_2(s)\big\rangle_{H}ds  \\
\ns\ds \q  + \mE \int_{t_{i+1}^n}^T
\big\langle Q^{m}(s)U(s,t_i^n) \xi_1,
v_2(s) \big\rangle_{H}ds.
\end{array}
\end{equation}
From \eqref{6.16eq1} and \eqref{6.16eq1xx},
it follows that
\begin{equation}\label{6.16eq1xxx}
\begin{array}{ll} \ds  \mE \big\langle
P^{m}(t_{i+1}^n)\xi_1,
x_2(t_{i+1}^n)\big\rangle_{H} -\mE
\int_{t_i^n}^{t_{i+1}^n} \big\langle
F_m(s) x_1(s), x_2(s) \big\rangle_{H}ds\\
\ns\ds = \mE \int_{t_i^n}^{t_{i+1}^n}
\big\langle P^{m}(s)K(s)x_1(s),
v_2(s)\big\rangle_{H}ds + \mE
\int_{t_i^n}^{t_{i+1}^n} \big\langle
Q^{m}(s)U(s,t_i^n) \xi_1, v_2(s)
\big\rangle_{H}ds,
\end{array}
\end{equation}
holds for any $i\in \{1,2,\cdots,n-1\}$,
$\xi_1\in L^{2q}_{\cF_{t^n_i}}(\O;H)$ and
$v_2\in
L^2_\dbF(t_i^n,t_{i+1}^n;L^{2q}(\O;H))$.

We choose a $\xi_1\in
L^{2q}_{\cF_{t_i^n}}(\O;H)$ with
$|\xi_1|_{L^{2q}_{\cF_{t_i^n}}(\O;H)}=1$
such that
$$
\big|Q^{m}(\cd)U(\cd,t_i^n)
\xi_1\big|_{L^{2}_\dbF(t_i^n,t_{i+1}^n;L^{\frac{2p}{p+1}}(\O;H))}
\geq
\frac{1}{2}\big|\!\big|Q^{m}(\cd)U(\cd,t_i^n)\big|\!\big|_{\cL(L^{2q}_{\cF_{t_i^n}}(\O;H),\;
L^{2}_\dbF(t_i^n,t_{i+1}^n;L^{\frac{2p}{p+1}}(\O;H))}.
$$
Next, we find a $v_2\in
L^2_\dbF(t_i^n,t_{i+1}^n;L^{2q}(\O;H))$
with
$|v_2|_{L^2_\dbF(t_i^n,t_{i+1}^n;L^{2q}(\O;H))}=1$
so that
$$
\mE \int_{t_i^n}^{t_{i+1}^n} \big\langle
Q^{m}(s)U(s,t_i^n)\xi_1, v_2(s)
\big\rangle_{H}ds \geq
\frac{1}{2}\big|Q^{n}(\cd)U(\cd,t_i^n)
\xi_1^n\big|_{L^{2}_\dbF(t_i^n,t_{i+1}^n;L^{\frac{2p}{p+1}}(\O;H))}.
$$
Hence,
\begin{equation}\label{6.16eq2}
\mE \int_{t_i^n}^{t_{i+1}^n} \big\langle
Q^{m}(s)U(s,t_i^n) \xi_1, v_2(s)
\big\rangle_{H}ds \geq
\frac{1}{4}\big|\!\big|Q^{m}(\cd)U(\cd,t_i^n)\big|\!\big|_{\cL(L^{2q}_{\cF_{t_i^n}}(\O;H),\;
L^{2}_\dbF(t_i^n,t_{i+1}^n;L^{\frac{2p}{p+1}}(\O;H))}.
\end{equation}
Also, it is easy to see that
\begin{equation}\label{6.16eq4}
\begin{array}{ll}\ds
\Big|P^{m}(t_{i+1}^n)\xi_1,
x_2(t_{i+1}^n)\big\rangle_{H} -\mE\!
\int_{t_i^n}^{t_{i+1}^n}\!\! \big\langle
F_m(s) x_1(s), x_2(s) \big\rangle_{H}ds -
\mE \!\int_{t_i^n}^{t_{i+1}^n}\!\!
\big\langle
P^m(s)K(s)x_1(s), v_2(s)\big\rangle_{H}ds \Big|\\
\ns\ds \leq\! C
\big(|P_T|_{L^p_{\cF_T}(\O;\cL(H))} +
|F|_{L^1_\dbF(0,T;L^p(\O;\cL(H)))}\big)\big(1
+
|(J,K)|_{(L^4_\dbF(0,T;L^\infty(\O;\cL(H))))^2}\big).
\end{array}
\end{equation}
Combining \eqref{6.16eq1}, \eqref{6.16eq2}
and \eqref{6.16eq4}, we find that
\begin{equation}\label{4.20-eq1}
\begin{array}{ll}\ds
\big|\!\big|Q^{m}(\cd)U(\cd,t_i^n)\big|\!\big|_{\cL(L^{2q}_{\cF_{t_i^n}}(\O;H),\;
L^{2}_\dbF(t_i^n,t_{i+1}^n;L^{\frac{2p}{p+1}}(\O;H)))}\\
\ns\ds \leq\!
C\big(|P_T|_{L^p_{\cF_T}(\O;\cL(H))} +
|F|_{L^1_\dbF(0,T;L^p(\O;\cL(H)))}\big)\big(1
+
|(J,K)|_{(L^4_\dbF(0,T;L^\infty(\O;\cL(H))))^2}\big).
\end{array}
\end{equation}
By \eqref{4.20-eq1} and Lemma \ref{cor1},
there exist a bounded, pointwisely defined
linear operator $ Q_{t_i^n}$ from
$L^{2q}_{\cF_{t_i^n}}(\O;H)$ to
$L^{2}_\dbF(t_i^n,t_{i+1}^n;L^{\frac{2p}{p+1}}(\O;H))$,
and a subsequence
$\{m_k^{(1)}\}_{k=1}^\infty$ of
$\{m\}_{m=1}^\infty$  such that
\begin{equation}\label{6.17eq1}
\mbox{(w)}\mbox{-}\lim_{k\to\infty}Q^{m_k^{(1)}}(\cd)U(\cd,t_i^n)\xi
= Q_{t_i^n}(\cd)\xi \q\mbox{ in }
L^{2}_\dbF(t_i^n,t_{i+1}^n;L^{\frac{2p}{p+1}}(\O;H)),\qq
\forall\; \xi\in L^4_{\cF_{t_i^n}}(\O;H).
\end{equation}

Since $ Q_{t_{i}^n}$ is pointwisely
defined, for $\ae (t,\o)\in
(t_{i}^n,t_{i+1}^n)\times\O$, there is a $
q_{t_{i}^n}(t,\o)\in \cL(H)$ such that
$$
\big(Q_{t_{i}^n} \xi\big)(t,\o)=
q_{t_{i}^n}(t,\o)\xi(\o),\qq \forall\;
\xi\in L^{2q}_{\cF_{t_{i}^n}}(\O;H).
$$
Let us define an operator $Q^{n}$ from
$\cH_n$ to
$L^2_\dbF(0,T;L^{\frac{2p}{p+1}}(\O;H))$ as
follows:
$$
\big(Q^{n}v\big)(t,\o)=\sum_{i=0}^{n-1}
\chi_{[t_{i}^n,t_{i+1}^n)}(t) q_{t_i}
(t,\o)h_i, \ \ae (t,\o)\in (0,T)\times \O,
$$
where $ v=\sum_{i=0}^{n-1}
\chi_{[t_{i}^n,t_{i+1}^n)}(\cd)U(\cd,t_{i}^n)h_i\in\cH_n$
with $h_i\in L^{2q}_{\cF_{t_{i}^n}}(\O;H)$.
It is easy to check that $ Q^{n} v\in
L^2_\dbF(0,T;L^{\frac{2p}{p+1}}(\O;H))$,
$Q^{n}$ is a pointwisely defined linear
operator from $\cH_n$ to
$L^2_\dbF(0,T;L^{\frac{2p}{p+1}}(\O;H))$
and
$$
\big|Q^{n}v\big|_{L^2_\dbF(0,T;L^{\frac{2p}{p+1}}(\O;H))}\leq
C|\tilde v|_{L^2_\dbF(0,T;L^{2q}(\O;H))}
$$
for $ \tilde v=\sum_{i=0}^{n-1}
\chi_{[t_i,t_{i+1}]}(\cd)h_i\in\cH_{n}$,
where $C$ is independent of $n$.

Also, for the above $v$, we have $
Q^{{m_{k}^{(1)}}}(s)v(s) =\sum_{i=0}^{n-1}
\chi_{[t_{i}^n,t_{i+1}^n)}Q^{{m_{k}^{(1)}}}(s)U(s,t_{i}^n)h_i$.
Hence,
$$
\ds Q^{{m_{k}^{(1)}}}(\cd)v(\cd)-
\big(Q^{n} v\big)(\cd)= \sum_{i=0}^{n-1}
\chi_{[t_i^n,t_{i+1}^n)}(\cd)\[
Q^{{m_{k}^{(1)}}}(\cd)U(s,t_i^n)h_i-\big(
Q_{t_i^n} h_i\big)(\cd)\].
$$
This gives that
\begin{equation}\label{s6eq1}
\mbox{(w)-}\lim_{k\to\infty}Q^{{m_{k}^{(1)}}}(\cd)v(\cd)=
Q^{n} v \q\mbox{ in }
L^{2}_\dbF(0,T;L^{\frac{2p}{p+1}}(\O;H)),\q
\forall\; v\in \cH_n.
\end{equation}

Similarly, one can find a subsequence
$\{m_{k}^{(2)}\}_{k=1}^\infty\subset\{m_{k}^{(1)}\}_{k=1}^\infty$
and a pointwisely defined linear operator
$\widehat Q^{n}$ from $\cH_n$ to
$L^2_\dbF(0,T;L^{\frac{2p}{p+1}}(\O;H))$
such that
\begin{equation}\label{s6eq1x}
\mbox{(w)-}\lim_{k\to\infty}Q^{{m_{k}^{(2)}}}(\cd)^*v(\cd)=
\widehat Q^n v \q\mbox{ in }
L^{2}_\dbF(0,T;L^{\frac{2p}{p+1}}(\O;H)),\q
\forall\; v\in \cH_n.
\end{equation}

For any $\xi_1,\xi_2\in
L^{2q}_{\cF_{0}}(\O;H)$, $u_1(\cd),
u_2(\cd)\in L^{2q}_{\dbF}(\O;L^2(0,T;H))$
and $v_1(\cd),v_2(\cd)\in \cH_n$, by
\eqref{s6eq1}--\eqref{s6eq1x}, it is easy
to see that
\begin{equation}\label{6.17eq2}
\begin{array}{ll}\ds
\lim_{k\to\infty}  \mE \int_{0}^T
\[\big\langle Q^{{m_{k}^{(2)}}}(s)v_1(s),
x_2(s) \big\rangle_{H}  +\big\langle
Q^{{m_{k}^{(2)}}}(s)x_1(s), v_2(s)
\big\rangle_{H}\]ds \\
\ns\ds  =  \mE \int_{0}^T
 \[\big\langle
\big(Q^{n} v_1\big)(s), x_2 (s)
\big\rangle_{H}+\big\langle x_1 (s),
\big(\widehat Q^{n} v_2\big)(s)
\big\rangle_{H}\]ds.
\end{array}
\end{equation}
On the other hand, from the proof of
\cite[Theorem 6.1]{LZ1}, one can show that
there exists a subsequence
$\{m_k^{(3)}\}_{k=1}^\infty$ of
$\{m_k^{(2)}\}_{m=1}^\infty$  such that
\begin{equation}\label{6.17euq2}
\begin{array}{ll}\ds
\lim_{k\to\infty}  \mE \int_{0}^T
\[\big\langle Q^{{m_{k}^{(3)}}}(s)v_1(s),
x_2(s) \big\rangle_{H}  +\big\langle
Q^{{m_{k}^{(3)}}}(s)x_1(s), v_2(s)
\big\rangle_{H}\]ds \\
\ns\ds  =  \mE \int_{0}^T \big\langle
v_1(s), \widehat Q^{(0)}(\xi_2,u_2,v_2)(s)
\big\rangle_{H}ds + \mE \int_0^T
\big\langle
 Q^{(0)}(\xi_1,u_1,v_1) (s),
v_2(s) \big\rangle_{H}ds .
\end{array}
\end{equation}
Combining \eqref{6.17eq2} and
\eqref{6.17euq2}, we obtain
\eqref{10.9eq2}. This completes the proof
of Theorem \ref{10.1th}.
\endpf


\section{Pontryagin-type maximum principle for controlled stochastic evolution equations}
\label{max}

In this section, for simplicity of the
presentation, we only consider the case
that $\{w(t)\}_{t\geq 0}$ is a standard one
dimensional Brownian motion.

Let $U$ be a separable metric space with
metric $d(\cd,\cd)$. Put
$$\cU[0,T] \triangleq \Big\{u(\cdot):\,
[0,T]\to U\;\Big|\; u(\cd) \mbox{ is
$\dbF$-adapted} \Big\}.$$ Consider the
following controlled (forward) stochastic
evolution equation:
\begin{eqnarray}\label{fsystem2}
\left\{
\begin{array}{lll}\ds
dx(t) = \big[Ax(t) +a(t,x(t),u(t))\big]dt + b(t,x(t),u(t))dw(t) &\mbox{ in }(0,T],\\
\ns\ds x(0)=x_0,
\end{array}
\right.
\end{eqnarray}
where $u(\cdot)\in \cU[0,T]$ and $x_0\in
L^8_{\cF_0}(\O;H)$. Similar to
\eqref{fsystem1}, $x(\cd)\equiv
x(\cd\,;x_0,u(\cdot))\!\in\!
C_\dbF([0,T];L^8(\O;H))$ is understood as a
mild solution to the equation
\eqref{fsystem2}.

Similar to \cite{LZ1}, we assume the
following three conditions:

\ms

\no{\bf (S1)} {\it Suppose that
$a(\cd,\cd,\cd),\,b(\cd,\cd,\cd):[0,T]\times
H\times U\to H$  are two maps such that for
$\f(t,x,u)=a(t,x,u),b(t,x,u)$, it holds
that: i) For any $(x,u)\in H\times U$, the
map $\f(\cd,x,u):[0,T]\to H$ is Lebesgue
measurable; ii) For any $(t,x)\in
[0,T]\times H$, the map $\f(t,x,\cd):U\to
H$ is continuous, and
\begin{equation}\label{ab0}
\left\{
\begin{array}{ll}\ds
|\f(t,x_1,u) - \f(t,x_2,u)|_H  \leq
C_L|x_1-x_2|_H,\q\forall\, (t,x_1,x_2,u)\in
[0,T]\times H\times
H\times U,\\
\ns\ds |\f(t,0,u)|_H \leq C_L, \q \forall\,
(t,u)\in [0,T]\times H\times H\times U;
\end{array}
\right.
\end{equation}}

\no{\bf (S2)} {\it Suppose that
$g(\cd,\cd,\cd):[0,T]\times H\times U\to
\dbR$ and $h(\cd):H\to \dbR$ are two
functionals such that for
$\psi(t,x,u)=g(t,x,u),h(x)$, it holds that:
i) For any $(x,u)\in H\times U$, the
function $\psi(\cd,x,u):[0,T]\to \dbR$ is
Lebesgue measurable; ii) For any $(t,x)\in
[0,T]\times H$, the function
$\psi(t,x,\cd):U\to \dbR$ is continuous,
and
\begin{equation}\label{gh} \left\{
\begin{array}{ll}\ds
|\psi(t,x_1,u) - \psi(t,x_2,u)|_{H}
 \leq C_L|x_1-x_2|_H,\q \forall\;
(t,x_1,x_2,u)\in [0,T]\times H\times H\times U,\\
\ns\ds |\psi(t,0,u)|_H  \leq C_L,\q
\forall\; (t,u)\in [0,T] \times U;
\end{array}
\right.
\end{equation}}

\no{\bf (S3)} {\it The map  $a(t,x,u)$ and
$b(t,x,u)$, and the functional $g(t,x,u)$
and $h(x)$ are $C^2$ with respect to $x$,
such that for
$\f(t,x,u)=a(t,x,u),b(t,x,u)$,
$\psi(t,x,u)= g(t,x,u),h(x)$, it holds that
 $\f_x(t,x,u)$, $\psi_x(t,x,u)$,
$\f_{xx}(t,x,u)$ and $\psi_{xx}(t,x,u)$ are
continuous with respect to $u$. Moreover,
\begin{equation}\label{ab1}
\left\{
\begin{array}{ll}\ds
|\!|\f_x(t,x,u)|\!|_{\cL(H)}+
|\psi_x(t,x,u) |_H \leq C_L,\q \forall\;
(t,x,u)\in
[0,T]\times H\times U,\\
\ns\ds |\!|\f_{xx}(t,x,u)|\!|_{\cL(H\times
H,\;H)} +|\!|\psi_{xx}(t,x,u)|\!|_{\cL(H)}
 \leq C_L,\q \forall\; (t,x,u)\in
[0,T]\times H\times U.
\end{array}
\right.
\end{equation}}

Define a cost functional $\cJ(\cdot)$ (for
the controlled equation \eqref{fsystem2})
as follows:
 \bel{jk1}
\cJ\big(u(\cdot)\big)\triangleq
\dbE\Big[\int_0^T g\big(t,x(t),u(t)\big)dt
+ h\big(x(T)\big)\Big],\qq\forall\;
u(\cdot)\in \cU[0,T].
 \ee
Let us consider the following optimal
control problem for \eqref{fsystem2}:

\ms

\no {\bf Problem (OP)} \q {\it Find a $\bar
u(\cdot)\in \cU[0,T]$ such that
 \bel{jk2}
\cJ \big(\bar u(\cdot)\big) =
\inf_{u(\cdot)\in \cU[0,T]} \cJ
\big(u(\cdot)\big).
 \ee
Any $\bar u(\cdot)$ satisfying (\ref{jk2})
is called an {\it optimal control}. The
corresponding state process $\bar x(\cdot)$
is called an {\it optimal state process}.
$\big(\bar x(\cdot),\bar u(\cdot)\big)$ is
called an {\it optimal pair}.}

\ms

There exist some works addressing the
Pontryagin-type maximum principle for
optimal controls of infinite dimensional
stochastic evolution equations (e.g.
\cite{Al-H1, Bensoussan2, HP1, TL, Zhou}
and the references therein).  However, most
of the previous works in this respect
addressed only to the case that either the
diffusion term does NOT depend on the
control variable (i.e., the map $b(t,x,u)$
in \eqref{fsystem2} is independent of $u$)
or the control region $U$ is convex.
Recently, this restriction was relaxed in
\cite{Du,Fuhrman,LZ1}. In both \cite{Du}
and \cite{Fuhrman}, the filtration $\dbF$
is assumed to be the natural one (generated
by the Brownian motion
$\{w(t)\}_{t\in[0,T]}$ and augmented by all
of the $\dbP$-null sets). Also, in
\cite{Du}, the authors assume that $A$ is a
strictly monotone operator; while in
\cite{Fuhrman}, the authors assume that
$H=L^2(D,\cD, \mu)$ (for a measure space
$(D,\cD,\mu)$ with finite measure $\mu$),
and the restriction of $\{S(t)\}_{t\geq 0}$
to the space $L^4(D,\cD,\mu)$ is a strongly
continuous analytic semigroup and the
domain of its infinitesimal generator is
compactly embedded in $L^4(D,\cD,\mu)$. On
the other hand, in \cite[Theorem 9.1]{LZ1},
a technical assumption $b_x(\cd,\bar
x(\cd),\bar u(\cd))\in
L^4_\dbF(0,T;L^\infty(\O;\cL(D(A))))$ is
imposed. The purpose of this section is to
establish a Pontryagin-type maximum
principle without any of the above
mentioned assumptions.

Define a function $\dbH:\; [0,T]\times H
\times U\times H\times H\to\dbR$ as
follows:
\begin{equation}\label{H}
\begin{array}{ll}\ds
\dbH(t,x,u,k_1,k_2) \= \big\langle k_1,a(t,x,u)  \big\rangle_H + \big\langle k_2, b(t,x,u)  \big\rangle_{H} - g(t,x,u),\\
\ns\ds \hspace{4cm} (t,x,u,k_1,k_2)\in
[0,T]\times H \times U\times H\times H.
\end{array}
\end{equation}
We have the following result.

\begin{theorem}\label{maximum p2}
Suppose that $H$ is a separable Hilbert
space, $L^p_{\cF_T}(\O;\dbC)$ ($1\leq p <
\infty$) is a separable Banach space, and
$U$ is a separable metric space. Let the
conditions (S1), (S2) and (S3) hold, and
$(\bar x(\cd),\bar u(\cd))$ be an optimal
pair for Problem (OP). Let
$\big(y(\cdot),Y(\cdot)\big)$ be the
transposition solution to the equation
\eqref{bsystem1} with $p=2$, and $y_T$ and
$f(\cd,\cd,\cd)$ given by
\begin{equation}\label{MP2-eq0}
\left\{
\begin{array}{ll} \ds
y_T = -h_x\big(\bar x(T)\big),\\\ns
 \ds f(t,y_1,y_2)=-a_x(t,\bar x(t),\bar
u(t))^*y_1 - b_x\big(t,\bar x(t),\bar
u(t)\big)^*y_2 + g_x\big(t,\bar x(t),\bar
u(t)\big). \end{array} \right.
\end{equation}
Assume that  $(P(\cd),Q^{(\cd)},\widehat
Q^{(\cd)})$ is the relaxed transposition
solution to the equation  \eqref{bsystem2}
in which $P_T$, $J(\cd)$, $K(\cd)$ and
$F(\cd)$ are given by
\begin{equation}\label{MP2-eq9}
\left\{
\begin{array}{ll} \ds P_T =
-h_{xx}\big(\bar x(T)\big),\q &\ds J(t) =
a_x(t,\bar x(t),\bar
u(t)), \\
\ns \ds K(t) =b_x(t,\bar x(t),\bar u(t)),
\q &\ds F(t)= -\dbH_{xx}\big(t,\bar
x(t),\bar u(t),y(t),Y(t)\big).
\end{array}
\right.
\end{equation}
Then,
\begin{equation}\label{MP2-eq1}
\begin{array}{ll}\ds
\Re\dbH\big(t,\bar x(t),\bar u(t),y(t),Y(t)\big) - \Re\dbH\big(t,\bar x(t),u,y(t),Y(t)\big) \\
\ns\ds \q - \frac{1}{2}\big\langle
P(t)\big[ b\big(t,\bar x(t),\bar
u(t)\big)-b\big(t,\bar x(t),u\big)
\big], b\big(t,\bar x(t),\bar u(t)\big)-b\big(t,\bar x(t),u\big) \big\rangle_H \\
\ns\ds\geq 0,\qq\q \,\ae [0,T]\times \O,\ \
\forall \,u \in U.
\end{array}
\end{equation}
\end{theorem}

{\it Proof}\,: We divide the proof into two
steps.

\ms

{\bf Step 1}. For each $\e>0$, let
$E_\e\subset [0,T]$ be a measurable set
with measure $\e$. Put
 $$
u^\e(t) = \left\{
\begin{array}{ll}
\ds \bar u(t), & t\in [0,T]\setminus E_\e,\\
\ns\ds u(t), & t\in E_\e,
\end{array}
\right.
 $$
where  $u(\cdot)$ is an arbitrary given
element in $\cU[0,T]$. Write
\begin{equation}\label{s7tatb1}
\left\{
\begin{array}{ll}
\ds a_1(t) = a_x(t,\bar x(t),\bar u(t)),
\qq b_1(t) = b_x(t,\bar
x(t),\bar u(t)),\qq g_1(t) = g_x(t,\bar x(t),\bar u(t)),\\
\ns\ds a_{11}(t) = a_{xx}(t,\bar x(t),\bar
u(t)), \qq b_{11}(t) = b_{xx}(t,\bar
x(t),\bar u(t)), \qq g_{11}(t)
= g_{xx}(t,\bar x(t),\bar u(t)),\\
\ns \ds \d a(t)  = a(t,\bar x(t), u(t)) -
a(t,\bar x(t),\bar u(t)),\qq\d b(t)  =
b(t,\bar x(t), u(t)) - b(t,\bar x(t),\bar
u(t)),
\\
\ns\ds\d g(t)  = g(t,\bar x(t), u(t)) -
g(t,\bar x(t),\bar u(t)), \qq \d b_1(t) =
b_x(t,\bar x(t), u(t)) - b_x(t,\bar
x(t),\bar u(t)).
\end{array}
\right.
\end{equation}

We introduce the following two stochastic
evolution equations:
\begin{equation}\label{MP2-equ2}
\left\{
\begin{array}{lll}\ds
dx_2^{\e } = \big[Ax_2^{\e } +
a_1(t)x_2^{\e } \big]dt + \big[ b_1(t)
x_2^{\e } + \chi_{E_\e} (t)\d b(t) \big]
dw(t) &\mbox{ in
}(0,T],\\
\ns\ds x_2^{\e }(0)=0
\end{array}
\right.
\end{equation}
and\footnote{Recall that, for any
$C^2$-function $f(\cd)$ defined on a Banach
space $X$ and $x_0\in X$, $f_{xx}(x_0)\in
\cL (X\times X,X)$. This means that, for
any $x_1,x_2\in X$,
$f_{xx}(x_0)(x_1,x_2)\in X$. Hence, by
\eqref{s7tatb1},
$a_{11}(t)\big(x_2^\e,x_2^\e\big)$ (in
\eqref{MP2-equ3}) stands for $a_{xx}(t,\bar
x(t),\bar
u(t))\big(x_2^\e(t),x_2^\e(t)\big)$. One
has a similar meaning for
$b_{11}(t)\big(x_2^\e,x_2^\e\big)$ and so
on.}
\begin{equation}\label{MP2-equ3}
\left\{
\begin{array}{lll}\ds
dx_3^{\e } = \[Ax_3^{\e } + a_1(t)x_3^{\e }
+ \chi_{E_\e}(t)\d a(t) +
\frac{1}{2}a_{11}(t)\big(x_2^{\e },x_2^{\e }\big) \]dt \\
\ns\ds \hspace{1cm} + \[ b_1(t) x_3^{\e } +
\chi_{E_\e} (t)\d b_1(t)x_2^{\e } +
\frac{1}{2}b_{11}(t)\big(x_2^{\e },x_2^{\e
}\big)\] dw(t) &\mbox{ in
}(0,T],\\
\ns\ds x_3^{\e }(0)=0.
\end{array}
\right.
\end{equation}
When $\e\to0$, by the proof of
\cite[Theorem 9.1]{LZ1}, we have
\begin{equation}\label{MP2-eq2}
\left\{
\begin{array}{ll}\ds
\ds |x_2^{\e }
|_{L^\infty_\dbF(0,T;L^{2}(\O;H))}=O(\sqrt{\e}),\\
\ns\ds |x_3^{\e }
|_{L^\infty_\dbF(0,T;L^{2}(\O;H))}=O(\e),
\end{array}
\right.
\end{equation}
and
\begin{equation}\label{s7eq5}
\begin{array}{ll}
\ds\cJ(u^\e(\cd)) -  \cJ(\bar u(\cd))\\
\ns\ds = \Re\mE \int_0^T\[ \big\langle
g_1(t),x_2^\e(t) + x_3^\e(t) \big\rangle_H
+ \frac{1}{2}\big\langle
g_{11}(t)x_2^\e(t),x_2^\e(t) \big\rangle_H
+ \chi_{E_\e}(t)\d g(t)
\] dt \\
\ns\ds\q + \Re\mE \big\langle h_x\big(\bar
x(T)\big), x_2^\e(T)+x_3^\e(T)
\big\rangle_H +
\frac{1}{2}\Re\mE\big\langle
h_{xx}\big(\bar
x(T)\big)x_2^\e(T),x_2^\e(T) \big\rangle_H
+ o(\e).
\end{array}
\end{equation}

By the definition of the transposition
solution to the equation \eqref{bsystem1}
(with $y_T$ and $f(\cd,\cd,\cd)$ given by
\eqref{MP2-eq0}), we obtain that
\begin{equation}\label{s7eq6}
-\mE\big\langle h_x(\bar
x(T))),x_2^\e(T)\big\rangle_H - \mE
\int_0^T \big\langle
g_1(t),x_2^\e(t)\big\rangle_H dt = \mE
\int_0^T\big\langle Y(t), \d
b(t)\big\rangle_H\chi_{E_\e}(t) dt
\end{equation}
and
\begin{equation}\label{s7eq7}
\begin{array}{ll}
\ds -\mE\big\langle h_x(\bar
x(T))),x_3^\e(T)\big\rangle_H - \mE
\int_0^T \big\langle
g_1(t),x_3^\e(t)\big\rangle_H dt
\\
\ns\ds = \mE \int_0^T \Big\{  \frac{1}{2}\[
\big\langle y(t),a_{11}(t)\big(x_2^\e(t),
x_2^\e(t)\big) \big\rangle_H + \big\langle
Y(t), b_{11}(t)\big(
x_2^\e(t), x_2^\e(t)\big) \big\rangle_H \] \\
\ns\ds \hspace{1.8cm} + \chi_{E_\e}(t)
\[ \big\langle y(t),\d a(t) \big\rangle_H +
\big\langle Y,\d b_1(t)x_2^\e(t)
\big\rangle_H
 \]\Big\}dt.
\end{array}
\end{equation}
According to
\eqref{MP2-eq2}--\eqref{s7eq7}, we conclude
that
\begin{equation}\label{s7eq8}
\begin{array}{ll} \ds
\cJ(u^\e(\cd)) -  \cJ(\bar u(\cd))\\
\ns\ds = \frac{1}{2}\Re\mE\int_0^T\[
\big\langle g_{11}(t)x_2^\e(t),
x_2^\e(t)\big\rangle_H - \big\langle
y(t),a_{11}(t)\big(x_2^\e(t),
x_2^\e(t)\big) \big\rangle_H \\
\ns\ds \q - \big\langle Y,
b_{11}(t)\big(x_2^\e(t),
x_2^\e(t)\big)\big\rangle_H
\]dt + \Re\mE\int_0^T \chi_{E_\e}(t)\[ \d g(t) - \big\langle
y(t),\d a(t)\big\rangle_H \\
\ns\ds \q-\big\langle Y(t),\d b(t)
\big\rangle_H
\]dt + \frac{1}{2}\Re\mE \big\langle h_{xx}\big(\bar
x(T)\big)x_2^\e(T), x_2^\e(T) \big\rangle_H
+ o(\e),\qq\hb{as }\e\to0.
\end{array}
\end{equation}

\ms

{\bf Step 2}. By the definition of the
relaxed transposition solution to the
equation \eqref{bsystem2} (with $P_T$,
$J(\cd)$, $K(\cd)$ and $F(\cd)$ given by
\eqref{MP2-eq9}), we obtain that
\begin{equation}\label{s7eq9}
\begin{array}{ll}\ds
-\mE\big\langle h_{xx}\big(\bar x(T)\big)
x_2^\e(T), x_2^\e(T) \big\rangle_H +
\mE\int_0^T \big\langle
\dbH_{xx}\big(t,\bar
x(t),\bar u(t),y(t),Y(t)\big) x_2^\e(t), x_2^\e(t) \big\rangle_H dt\\
\ns\ds = \dbE\int_0^T
\chi_{E_\e}(t)\big\langle b_1(t)x_2^\e(t),
P(t)^*\d b(t)\big\rangle_{H} dt +
\dbE\int_0^T \chi_{E_\e}(t)\big\langle
P(t)\d b(t),
b_1(t)x_2^\e(t)\big\rangle_{H} dt\\
\ns\ds \q  + \dbE\int_0^T
\chi_{E_\e}(t)\big\langle P(t)\d b(t), \d
b(t)\big\rangle_{H} dt + \dbE\int_0^T
\chi_{E_\e}(t)\big\langle
\d b(t),\widehat Q^{(0)}(0,0,\chi_{E_\e}\d b)(t)\big\rangle_{H} dt \\
\ns\ds \q + \dbE\int_0^T
\chi_{E_\e}(t)\big\langle Q^{(0)}(0,0,\d
b)(t),\d b(t)\big\rangle_{H} dt.
 \ea
\end{equation}

Now, we estimate the terms in the right
hand side of \eqref{s7eq9}. By
\eqref{MP2-eq2}, we have
\begin{equation}\label{s7eq9.1}
\begin{array}{ll}\ds
\Big|\dbE\int_0^T \chi_{E_\e}(t)\big\langle
b_1(t)x_2^\e(t), P(t)^*\d
b(t)\big\rangle_{H} dt+\dbE\int_0^T
\chi_{E_\e}(t)\big\langle P(t)\d b(t),
b_1(t)x_2^\e(t)\big\rangle_{H} dt\Big|=
o(\e).
\end{array}
\end{equation}

In what follows, for any $\tau\in [0,T)$,
we choose $E_{\e}=[\tau,\tau+\e]\subset
[0,T]$.

By Proposition \ref{5.13-prop1}, we can
find a sequence $\{\beta_n\}_{n=1}^\infty$
such that $\b_n\in\cH_n$ (Recall \eqref{cH}
for the definition of $\cH_n$) and $
\lim_{n\to\infty}\beta_n = \d b $ in $
L^2_\dbF(0,T;L^4(\O;H))$. Hence, for some
positive constant $C(x_0)$ (depending on
$x_0$),
\begin{equation}\label{10.9qqq3}
|\beta_n|_{L^2_\dbF(0,T;L^4(\O;H))}\le
C(x_0)<\infty,\qq\forall\;n\in\dbN,
\end{equation}
and there is a subsequence
$\{n_k\}_{k=1}^\infty\subset \{
n\}_{n=1}^\infty$ such that
\begin{equation}\label{s7eq9.2-11}
\lim_{k\to\infty} |\b_ {n_k}(t)-\d
b(t)|_{L^4_{\cF_t}(\O;H)} = 0\q\mbox{ for
}\ae t\in [0,T].
\end{equation}

Denote by $Q^{n_k}$ and $\widehat Q^{n_k}$
the corresponding pointwisely defined
linear operators from $\cH_{n_k}$ to
$L^2_\dbF(0,T;L^{\frac{4}{3}}(\O;H))$,
given in Theorem \ref{10.1th}.

Consider the following equation:
\begin{equation}\label{s7fsystem3.1x}
\left\{
\begin{array}{lll}\ds
dx_{2,n_k}^{\e} = \big[Ax_{2,n_k}^{\e} +
a_1(t)x_{2,n_k}^{\e} \big]dt + \big[ b_1(t)
x_{2,n_k}^{\e} + \chi_{E_{\e}} (t)\b_
{n_k}(t) \big] dw(t) &\mbox{ in
}(0,T],\\
\ns\ds x_{2,n_k}^{\e}(0)=0.
\end{array}
\right.
\end{equation}
We have
\begin{equation}\label{s7th max eq1x}
\begin{array}{ll}\ds
\mE|x_{2,n_k}^{\e}(t)|^4_H \\
\ns\ds = \mE\Big| \int_0^t
S(t-s)a_1(s)x_{2,n_k}^{\e}(s) ds + \int_0^t
S(t-s)b_1(s)x_{2,n_k}^{\e}(s)
dw(s) \\
\ns\ds \qq +
\int_0^t S(t-s)\chi_{E_{\e}}(s)\b_ {n_k}(s) dw(s)\Big|_H^4\\
\ns\ds  \leq C \bigg[\mE\Big| \int_0^t
S(t-s)a_1(s)x_{2,n_k}^{\e}(s) ds \Big|_H^4
+ \mE\Big|\int_0^t S(t-s)b_1(s)x_{2,n_k}^{\e}(s) dw(s) \Big|_H^4 \\
\ns\ds \q  + \mE\Big| \int_0^t S(t-s)\chi_{E_{\e}}(s)\b_ {n_k}(s) dw(s)\Big|_H^4\bigg] \\
\ns\ds   \leq C\[\int_0^t
\mE|x_{2,n_k}^{\e}(s)|_H^4 ds  + \e
\int_{E_{\e}}\mE|\beta_{n_k}(s)|_{H}^4ds
\].
\end{array}
\end{equation}
By \eqref{10.9qqq3} and thanks to
Gronwall's inequality, \eqref{s7th max
eq1x} leads to
\begin{equation}\label{s7th max eq2x}
|x_{2,n_k}^{\e}(\cd)|^4_{L^\infty_\dbF(0,T;L^4(\O;H))}
\leq C(x_0,k)\e^2.
\end{equation}
Here and henceforth, $C(x_0,k)$ is a
generic constant (depending on $x_0$,
 $k$, $T$, $A$ and $C_L$), which may
be different from line to line. For any
fixed $ k\in\dbN$, since $Q^{n_k}\b_
{n_k}\in
L^2_{\dbF}(0,T;L^{\frac{4}{3}}(\O;H))$, by
\eqref{s7th max eq2x}, we find that
\begin{equation}\label{s7eq9.3}
\begin{array}{ll}\ds
\Big|\dbE\int_0^T
\chi_{E_{\e}}(t)\big\langle \big(Q^{n_k}
\b_
{n_k}\big)(t),x_{2,n_k}^{\e}(t)\big\rangle_{H}
dt \Big|\\
\ns\ds \leq
|x_{2,n_k}^{\e}(\cd)|_{L^\infty_\dbF(0,T;L^4(\O;H))}
\int_{E_{\e}}\big|\big(Q^{n_k} \b_ {n_k}\big)(t)\big|_{L^{\frac{4}{3}}_{\cF_t}(\O;H)}dt\\
\ns\ds \leq C(x_0,k)\sqrt{{\e}}
\int_{E_{\e}}\big|\big(Q^{n_k} \b_
{n_k}\big)(t)\big|_{L^{\frac{4}{3}}_{\cF_t}(\O;H)}dt\\
\ns\ds = o({\e}), \qq\hbox{as }\e\to0.
\end{array}
\end{equation}
Similarly,
\begin{equation}\label{s7eq9.3x}
\Big|\dbE\int_0^T
\chi_{E_{\e}}(t)\big\langle
x_{2,n_k}^{\e}(t),\big(\widehat Q^{n_k}\b_
{n_k}\big)(t)\big\rangle_{H} dt \Big| =
o({\e}), \qq\hbox{as }\e\to0.
\end{equation}

From \eqref{10.9eq2} in Theorem
\ref{10.1th}, and noting that both
$Q^{n_k}$ and $\widehat Q^{n_k}$ are
pointwisely defined, we arrive at the
following equality:
\begin{equation}\label{wws1}
\begin{array}{ll}
\ds \mE \int_{0}^T \big\langle
\chi_{E_{\e}}(t)\b_ {n_k}(t), \widehat
Q^{(0)}(0,0,\chi_{E_{\e}} \b_{n_k})(t)
\big\rangle_{H}dt +  \mE \int_{0}^T
\big\langle
Q^{(0)}(0,0,\chi_{E_{\e}} \b_ {n_k}) (t), \chi_{E_{\e}}\b_ {n_k}(t) \big\rangle_{H}dt \\
\ns\ds =\mE \int_{0}^T
\chi_{E_{\e}}\[\big\langle \big(Q^{n_k} \b_
{n_k}\big)(t), x_{2,n_k}^{\e} (t)
\big\rangle_{H}+\big\langle x_{2,n_k}^{\e}
(t), \big(\widehat Q^{n_k} \b_
{n_k}\big)(t) \big\rangle_{H}\]dt.
\end{array}
\end{equation}
Hence,
\begin{equation}\label{wws2}
\begin{array}{ll}\ds
\mE \int_{0}^T \big\langle
\chi_{E_{\e}}(t)\d b(t), \widehat
Q^{(0)}(0,0,\chi_{E_{\e}}\d b)(t)
\big\rangle_{H}dt + \mE \int_{0}^T
\big\langle
 Q^{(0)}(0,0,\chi_{E_{\e}}\d
b) (t), \chi_{E_{\e}}(t)\d
b(t) \big\rangle_{H}dt \\
\ns\ds \q-\mE \int_{0}^T \chi_{E_{\e}}(t)
\[\big\langle \big(Q^{n_k} \b_ {n_k}\big)(t),
x_{2,n_k}^{\e}(t)
\big\rangle_{H}+\big\langle x_{2,n_k}^{\e}
(t), \big(\widehat Q^{n_k} \b_
{n_k}\big)(t)
\big\rangle_{H}\]dt\\
\ns\ds = \mE \int_{0}^T \big\langle
\chi_{E_{\e}}(t)\d b(t), \widehat
Q^{(0)}(0,0,\chi_{E_{\e}} \d b)(t)
\big\rangle_{H}dt + \mE \int_{0}^T
\big\langle
 Q^{(0)}(0,0,\chi_{E_{\e}} \d
b) (t), \chi_{E_{\e}}(t)\d
b(t) \big\rangle_{H}dt \\
\ns\ds \q\!\! -\mE\!\int_{0}^T\!\big\langle
\chi_{E_{\e}}\!(t)\b_ {n_k}\!(t), \widehat
Q^{(0)}(0,0,\chi_{E_{\e}}\b_{n_k})(t)\big\rangle_{H}dt
\! - \!\mE\!\int_{0}^T\!\big\langle
Q^{(0)}(0,0,\chi_{E_{\e}} \b_{n_k} ) (t),
\chi_{E_{\e}}\!(t)\b_{n_k}\!(t)
\big\rangle_{H}dt.
\end{array}
\end{equation}
It is easy to see that
\begin{equation}\label{s7eq9.2xx}
\begin{array}{ll}\ds\Big|\mE \int_{0}^T
\big\langle \chi_{E_{\e}}(t)\d b(t),
\widehat Q^{(0)}(0,0,\chi_{E_{\e}} \d b)(t)
\big\rangle_{H}dt  - \mE \int_{0}^T
\big\langle \chi_{E_{\e}}(t)\b_ {n_k}(t),
\widehat Q^{(0)}(0,0,\chi_{E_{\e}} \b_
{n_k})(t)
\big\rangle_{H}dt \Big|\\
\ns\ds \leq \Big|\mE \int_{0}^T \big\langle
\chi_{E_{\e}}(t)\d b(t), \widehat
Q^{(0)}(0,0,\chi_{E_{\e}} \d b)(t)
\big\rangle_{H}dt - \mE \int_{0}^T
\big\langle \chi_{E_{\e}}(t)\d b(t),
\widehat Q^{(0)}(0,0,\chi_{E_{\e}} \b_
{n_k})(t)
\big\rangle_{H}dt \Big|\\
\ns\ds \q + \Big|\mE\! \int_{0}^T\!\!
\big\langle \chi_{E_{\e}}(t)\d b(t),
\widehat Q^{(0)}(0,0,\chi_{E_{\e}} \b_
{n_k})(t) \big\rangle_{H}dt\! - \!\mE\!
\int_{0}^T\!\! \big\langle
\chi_{E_{\e}}(t)\b_ {n_k}(t), \widehat
Q^{(0)}(0,0,\chi_{E_{\e}} \b_ {n_k})(t)
\big\rangle_{H}dt \Big|.
\end{array}
\end{equation}
From \eqref{s7eq9.2-11} and the density of
the Lebesgue points, we find that for
$\ae\tau\in [0,T)$, it holds that
\begin{equation}\label{s7eq9.3-1}
\begin{array}{ll}\ds
\lim_{k\to\infty}\lim_{\e\to
0}\frac{1}{{\e}}\Big|\mE \int_{0}^T
\big\langle \chi_{E_{\e}}(t)\d
b(t),\widehat Q^{(0)}(0,0,\chi_{E_{\e}} \d
b)(t) \big\rangle_{H}dt
\\
\ns\ds \qq\qq\q - \mE \int_{0}^T
\big\langle
\chi_{E_{\e}}(t)\d b(t),\widehat Q^{(0)}(0,0,\chi_{E_{\e}} \b_ {n_k})(t) \big\rangle_{H}dt \Big|\\
\ns\ds \leq \lim_{k\to\infty}\lim_{\e\to
0}\frac{1}{{\e}}\[\int_0^T \chi_{E_{\e}}(t)
\(\mE|\d b(t)|^4_{H}\)^{\frac{1}{2}} dt
\Big]^{\frac{1}{2}}|\widehat
Q^{(0)}(0,0,\chi_{E_{\e}} (\d b-\b_ {n_k}))
|_{L^2_\dbF(0,T;L^{\frac{4}{3}}(\O;H))}
\\
\ns\ds \leq C\lim_{k\to\infty}\lim_{\e\to
0}\frac{1}{{\e}}
\[\int_0^T \chi_{E_{\e}}(t) \(\mE|\d b(t) |^4_{H}\)^{\frac{1}{2}} dt
\Big]^{\frac{1}{2}}\big|\chi_{E_{\e}} (\d b-\b_ {n_k})\big|_{L^2_\dbF(0,T;L^4(\O;H))}\\
\ns\ds \leq C\lim_{k\to\infty}\lim_{\e\to
0}\frac{|\d
b(\tau)|_{L^4_{\cF_\tau}(\O;H)}}{\sqrt{{\e}}}\[\int_0^T
\chi_{E_{\e}}(t) \(\mE|\d b(t) - \b_
{n_k}(t)|^4_{H}\)^{\frac{1}{2}} dt
\Big]^{\frac{1}{2}}\\
\ns\ds = C\lim_{k\to\infty}\lim_{\e\to
0}|\d
b(\tau)|_{L^4_{\cF_\tau}(\O;H)}\[\frac{1}{{\e}}\int_\tau^{\tau+{\e}}
|\d b(t) - \b_
{n_k}(t)|_{L^4_{\cF_t}(\O;H)}^2 dt
\Big]^{\frac{1}{2}}\\
\ns\ds = C\lim_{k\to\infty}|\d
b(\tau)|_{L^4_{\cF_\tau}(\O;H)}|\d b(\tau)
- \b_ {n_k}(\tau)|_{L^4_{\cF_\tau}(\O;H)}
\\\ns\ds= 0.
\end{array}
\end{equation}
Similarly,
\begin{equation}\label{s7eq9.3-1x}
\begin{array}{ll}\ds
\lim_{k\to\infty}\lim_{\e\to
0}\frac{1}{{\e}}\Big|\mE \int_{0}^T
\big\langle \chi_{E_{\e}}(t)\d b(t),
\widehat Q^{(0)}(0,0,\chi_{E_{\e}} \b_
{n_k})(t)
\big\rangle_{H}dt \\
\ns\ds \qq\qq\q - \mE  \int_{0}^T
\big\langle \chi_{E_{\e}}(t)\b_ {n_k}(t),
\widehat Q^{(0)}(0,0,\chi_{E_{\e}} \b_
{n_k})(t)
\big\rangle_{H}dt \Big|\\
\ns\ds \leq \lim_{k\to\infty}\lim_{\e\to
0}\frac{1}{{\e}}\big|\widehat
Q^{(0)}(0,0,\chi_{E_{\e}} \b_
{n_k})\big|_{L^2_\dbF(0,T;L^{\frac{4}{3}}(\O;H))}
\[\int_0^T \chi_{E_{\e}}(t) \(\mE|\d b(t) - \b_ {n_k}(t)|^4_{H}\)^{\frac{1}{2}} dt
\Big]^{\frac{1}{2}}\\
\ns\ds \leq C\lim_{k\to\infty}\lim_{\e\to
0}\frac{1}{{\e}}\big|\chi_{E_{\e}} \b_
{n_k}\big|_{L^2_\dbF(0,T;L^4(\O;H))}
\[\int_0^T \chi_{E_{\e}}(t) \(\mE|\d b(t) - \b_ {n_k}(t)|^4_{H}\)^{\frac{1}{2}} dt
\Big]^{\frac{1}{2}}\\
\ns\ds \leq C\lim_{k\to\infty}\lim_{\e\to
0}\frac{1}{{\e}}\Big\{\big|\chi_{E_{\e}} \d
b\big|_{L^2_\dbF(0,T;L^4(\O;H))}
\[\int_0^T \chi_{E_{\e}}(t) \(\mE|\d b(t) - \b_ {n_k}(t)|^4_{H}\)^{\frac{1}{2}} dt
\Big]^{\frac{1}{2}}\\\ns\ds\qq\qq\qq\qq\q+\int_0^T \chi_{E_{\e}}(t) \(\mE|\d b(t) - \b_ {n_k}(t)|^4_{H}\)^{\frac{1}{2}} dt\Big\}\\
\ns\ds \leq C\lim_{k\to\infty}\lim_{\e\to
0}\Big\{\frac{|\d
b(\tau)|_{L^4_{\cF_\tau}(\O;H)}}{\sqrt{{\e}}}\[\int_0^T
\chi_{E_{\e}}(t) \(\mE|\d b(t) - \b_
{n_k}(t)|^4_{H}\)^{\frac{1}{2}} dt
\Big]^{\frac{1}{2}}\\\ns\ds\qq\qq\qq\q+\frac{1}{{\e}}\int_0^T
\chi_{E_{\e}}(t) \(\mE|\d b(t) - \b_
{n_k}(t)|^4_{H}\)^{\frac{1}{2}} dt
\Big\}\\
\ns\ds = C\lim_{k\to\infty}\lim_{\e\to
0}\Big\{|\d
b(\tau)|_{L^4_{\cF_\tau}(\O;H)}\[\frac{1}{{\e}}\int_\tau^{\tau+{\e}}
|\d b(t) - \b_
{n_k}(t)|_{L^4_{\cF_t}(\O;H)}^2 dt
\Big]^{\frac{1}{2}}\\\ns\ds\qq\qq\qq\q+\frac{1}{{\e}}\int_\tau^{\tau+{\e}}
|\d b(t) - \b_
{n_k}(t)|_{L^4_{\cF_t}(\O;H)}^2 dt
\Big\}\\
\ns\ds = C\lim_{k\to\infty}\big[|\d
b(\tau)|_{L^4_{\cF_\tau}(\O;H)}|\d b(\tau)
- \b_
{n_k}(\tau)|_{L^4_{\cF_\tau}(\O;H)}+|\d
b(\tau) - \b_
{n_k}(\tau)|_{L^4_{\cF_\tau}(\O;H)}^2\big]
\\\ns\ds= 0.
\end{array}
\end{equation}
From \eqref{s7eq9.2xx}--\eqref{s7eq9.3-1x},
we find that
\begin{equation}\label{s7eq9.3-1xx}
\begin{array}{ll}\ds
\lim_{k\to\infty} \lim_{\e\to
0}\frac{1}{{\e}}\Big|\mE \int_{0}^T
\big\langle \chi_{E_{\e}}(t)\d b(t),
\widehat Q^{(0)}(0,0,\chi_{E_{\e}} \d b)(t)
\big\rangle_{H}dt \\
\ns\ds \qq - \mE  \int_{0}^T \big\langle
\chi_{E_{\e}}(t)\b_ {n_k}(t), \widehat
Q^{(0)}(0,0,\chi_{E_{\e}} \b_ {n_k})(t)
\big\rangle_{H}dt \Big| \\
\ns\ds= 0.
\end{array}
\end{equation}
By a similar argument, we obtain that
\begin{equation}\label{s7eq9.3-1xxx}
\begin{array}{ll}\ds
\lim_{k\to\infty} \lim_{\e\to
0}\frac{1}{{\e}}\Big|\mE \int_{0}^T
\big\langle Q^{(0)}(0,0,\chi_{E_{\e}} \d
b)(t),\chi_{E_{\e}}(t)\d b(t)
\big\rangle_{H}dt \\
\ns\ds \qq - \mE  \int_{0}^T \big\langle
Q^{(0)}(0,0,\chi_{E_{\e}} \b_
{n_k})(t),\chi_{E_{\e}}(t)\b_ {n_k}(t)
\big\rangle_{H}dt \Big| \\
\ns\ds= 0.
\end{array}
\end{equation}

From \eqref{s7eq9.3}--\eqref{wws2} and
\eqref{s7eq9.3-1xx}--\eqref{s7eq9.3-1xxx},
we obtain that
\begin{equation}\label{s7eq9.2-111}
 \ba{ll}\ds
\Big|\dbE\int_0^T
\chi_{E_{\e}}(t)\big\langle \d
b(t),\widehat Q^{(0)}(0,0,\chi_{E_{\e}}\d
b)(t)\big\rangle_{H} dt + \dbE\int_0^T
\chi_{E_{\e}}(t)\big\langle Q^{(0)}(0,0,\d
b)(t),\d b(t)\big\rangle_{H} dt\Big|
\\\ns\ds=o({\e}), \qq\hbox{as }\e\to0.
 \ea
\end{equation}

Combining \eqref{s7eq8}, \eqref{s7eq9},
\eqref{s7eq9.1} and \eqref{s7eq9.2-111}, we
end up with
\begin{equation}\label{s7eq10}
 \ba{ll}\ds
\cJ(u^{\e}(\cd)) -  \cJ(\bar
u(\cd))\\\ns\ds = \Re\mE\int_0^T \[
 \d g(t) - \big\langle
y(t),\d a(t)\big\rangle_H -\big\langle
Y(t),\d b(t) \big\rangle_H -
\frac{1}{2}\big\langle P(t)\d b(t), \d b(t)
\big\rangle_H
 \]\chi_{E_{\e}}(t)dt + o({\e}).
 \ea
\end{equation}
Since $\bar u(\cd)$ is the optimal control,
$\cJ(u^{\e}(\cd)) - \cJ(\bar u(\cd))\geq
0$. Thus,
\begin{equation}\label{s7eq11}
\Re\mE\int_0^T \chi_{E_{\e}}(t)\[
  \big\langle
y(t),\d a(t)\big\rangle_H +\big\langle
Y(t),\d b(t) \big\rangle_H -\d g(t)+
\frac{1}{2}\big\langle P(t)\d b(t), \d b(t)
\big\rangle_H
 \]dt \leq o({\e}),
\end{equation}
as $\e\to0$.

Finally, by \eqref{s7eq11}, we obtain
\eqref{MP2-eq1}. This completes the proof
of Theorem \ref{maximum p2}.
\endpf

\section*{Acknowledgement}

This work was partially supported by the
National Basic Research Program of China
(973 Program) under grant 2011CB808002, the
NSF of China under grants 11101070,
11221101 and 11231007, the PCSIRT under
grant IRT1273 and the Chang Jiang Scholars
Program (from the Chinese Education
Ministry), and the Grant MTM2011-29306
(from the Spanish Science and Innovation
Ministry).


\begin{thebibliography}{99}


\bibitem{Al-H1}A. Al-Hussein. \it Backward stochastic partial differential equations
driven by infinite dimensional martingales
and applications. \sl Stochastics. \rm {\bf
81} (2009), 601--626.

\bibitem{Bensoussan2}A. Bensoussan. \it Stochastic maximum principle for distributed
parameter systems. \sl J. Franklin Inst.
\rm {\bf 315} (1983), 387--406.

\bibitem{Du} K.~Du and Q.~Meng. \it A maximum
principle for optimal control of stochastic
evolution equations. \sl SIAM J. Control
Optim. \rm{\bf 51} (2013),  4343--4362.

\bibitem{Fuhrman} M.~Fuhrman, Y.~Hu and
G.~Tessitore. \it Stochastic maximum
principle for optimal control of SPDEs. \sl
Appl. Math. Optim. \rm {\bf 68} (2013),
181--217.

\bibitem{HP1} Y. Hu and S. Peng. \it Maximum principle for semilinear stochastic evolution control
systems. \sl Stochastics Stochastics Rep.
\rm {\bf 33} (1990), 159--180.

\bibitem{HP2} Y. Hu and S. Peng. \it Adapted solution of backward semilinear stochastic evolution
equations. \sl Stoch. Anal. Appl. \rm {\bf
9} (1991), 445--459.

\bibitem{Lions} J.-L.~Lions. \sl Contr\^olabilit\'e exacte, perturbations et
stabilisation de syst\`emes distribu\'es,
Tome 1. \rm Recherches en Math\'ematiques
Appliqu\'ees, Tome~8. Masson, Paris, 1988.

\bibitem{Lions1} J.-L.~Lions and
E.~Magenes. \sl Non-homogeneous Boundary
Value Problems and Applications, Vol. I.
\rm Die Grundlehren der mathematischen
Wissenschaften, Band {\bf 181}.
Springer-Verlag, New York, 1972.

\bibitem{LYZ}Q. L\"{u},  J.~Yong and X.~Zhang. \it Representation of It\^o integrals by Lebesgue/Bochner
integrals. \sl J. Eur. Math. Soc. \rm {\bf
14} (2012), 1795--1823.

\bibitem{LZ} Q. L\"{u} and X.~Zhang. \it  Well-posedness of
backward stochastic differential equations
with general filtration. \sl J.
Differential Equations. \rm{\bf 254}
(2013), 3200--3227.

\bibitem{LZ1} Q. L\"{u} and X.~Zhang. \it  General Pontryagin-type stochastic
maximum principle and backward stochastic
evolution equations in infinite dimensions.
\rm Springer Briefs in Mathematics,
Springer, New York, 2014. (See also
http://arxiv.org/abs/1204.3275)


\bibitem{MY}J. Ma and J. Yong. \sl Forward-Backward Stochastic Differential
Equations and Their Applications. \rm
Lecture Notes in Math. vol. {\bf 1702}.
Springer-Verlag, New York, 1999.

\bibitem{MM}N. I. Mahmudova and M. A. McKibben. \it On backward stochastic evolution equations in Hilbert spaces and
optimal control. \sl Nonlinear Anal. \rm
{\bf 67} (2007), 1260--1274.

\bibitem{Peng1} S.~Peng. \it A general stochastic maximum principle for optimal
control problems. \sl SIAM J. Control
Optim. \rm {\bf 28} (1990), 966--979.

\bibitem{TL} S. Tang and X. Li. \it Maximum principle for optimal
control of distributed parameter stochastic
systems with random jumps. \rm In: \sl
Differential Equations, Dynamical Systems,
and Control Science. \rm Edited by K. D.
Elworthy,  W. N. Everitt and E. B. Lee.
Lecture Notes in Pure and Appl. Math. vol.
{\bf 152}. Dekker, New York, 1994,
867--890.

\bibitem{YZ} J.~Yong and X. Y.~Zhou. \sl Stochastic Controls: Hamiltonian Systems and
HJB Equations. \rm Springer-Verlag, New
York, 1999.

\bibitem{Zhou} X. Y.~Zhou. \it On the necessary conditions of optimal controls for
stochastic partial differential equations.
\sl SIAM J. Control Optim. \rm {\bf 31}
(1993), 1462--1478.

\end{thebibliography}
\end{document}